\documentclass[11pt,a4paper]{amsart} 
 
\usepackage[utf8]{inputenc}
\usepackage[a4paper]{geometry}
\usepackage{latexsym}
\usepackage{kpfonts} 
\usepackage{amsmath,amssymb,amsfonts,amscd,amsthm}
\usepackage{mathtools}
\usepackage{physics}
\usepackage[mathscr]{eucal}
\usepackage{array}
\usepackage{hyperref}

\usepackage[style = alphabetic, backend = bibtex]{biblatex}
\theoremstyle{plain}
\newtheorem{theorem}{Theorem}
\newtheorem{proposition}{Proposition}
\newtheorem{lemma}{Lemma}

\theoremstyle{definition}

\theoremstyle{remark}

\newtheorem{remark}[theorem]{Remark}

\DeclareMathOperator{\vol}{vol}
\DeclareMathOperator{\disc}{disc}
\DeclareMathOperator{\nr}{nr}
\DeclareMathOperator{\SL}{SL}

\numberwithin{equation}{section}
\numberwithin{proposition}{section}
\numberwithin{lemma}{section}

\title[Hybrid sup-norm bounds in higher rank]{Hybrid bounds for the sup-norm of automorphic forms in higher rank}
\author{Radu Toma}
\address{Mathematisches Institut, Endenicher Allee 60, 53115 Bonn, Germany}
\email{toma@math.uni-bonn.de}
\thanks{Author supported by Germany Excellence Strategy grant EXC-2047/1-390685813 through the Bonn International Graduate School of Mathematics}

\begin{document}

\begin{abstract}
Let $A$ be a central division algebra of prime degree $p$ over $\mathbb{Q}$. We obtain subconvex hybrid bounds, uniform in both the eigenvalue and the discriminant, for the sup-norm of Hecke-Maaß forms on the compact quotients of $\operatorname{SL}_p(\mathbb{R})/\operatorname{SO}(p)$ by unit groups of orders in $A$. The exponents in the bounds are explicit and polynomial in $p$. We also prove subconvex hybrid bounds in the case of certain Eichler-type orders in division algebras of arbitrary odd degree.
\end{abstract}

\maketitle

\section{Introduction}

\subsection{Motivation and historical context}

The sup-norm problem arises as a natural question in analysis and quantum physics and has received considerable attention in the number theory community. It is the problem of bounding the $L^\infty$-norm of eigenfunctions on Riemannian manifolds in terms of their $L^2$-norm. To make this a reasonable endeavour, one chooses some parameters for the eigenfunctions and estimates the quotient of the two norms while these parameters vary.

For example, let $X$ be a compact Riemannian manifold of dimension $n$ and let $\phi$ be an $L^2$-normalised eigenfunction of the Laplacian $\Delta_X$ with eigenvalue $\lambda > 0$. Motivated by semi-classical analysis, one would like to bound $\norm{\phi}_\infty$ in terms of $\lambda$ when $\lambda \rightarrow \infty$. In general, local analysis gives the sharp baseline bound
\begin{displaymath}
\norm{\phi}_\infty \ll \lambda^{(n-1)/4 + \varepsilon},
\end{displaymath}
for large enough $\lambda$. The bound is attained on the round $n$-spheres.

If $\phi$ is assumed to be an eigenfunction for a larger algebra of operators, then we can expect better bounds. Indeed, if $X = \Gamma \backslash S$ is a compact locally symmetric space of rank $r$ and $G(S)$ is the groups of isometries of the symmetric space $S$, then one can consider the algebra of $G(S)$-invariant differential operators. This algebra is generated by $r$ operators, including the Laplacian. If $\phi$ is an $L^2$-normalised joint eigenfunction of these operators, then\footnote{Here and in the rest of this article, the implied constants are allowed to depend on $\varepsilon$.}
\begin{equation}\label{eq-local}
\norm{\phi}_\infty \ll \lambda^{(n-r)/4 + \varepsilon}.
\end{equation}
For more details on the above paragraphs, see \cite{sarnak}.

In special cases when $X$ has arithmetic structure, we expect to obtain better, so-called \emph{subconvex} bounds, i.e.\ an exponent $(n-r)/4 - \delta$ in \eqref{eq-local} with $\delta > 0$. In these cases, there is an additional algebra of commutative normal Hecke operators, which commute with the differential operators above. In this arithmetic setting, the sup-norm problem is to find a subconvex bound for $\norm{\phi}_\infty$, where $\phi$ is a joint eigenfunction of the invariant differential operators and the Hecke algebra. The prototype of such a result is due to Iwaniec and Sarnak \cite{IS} in the case of $X = \Gamma \backslash \mathfrak{h}^2$, where $\mathfrak{h}^2$ is the hyperbolic plane and $\Gamma \leq \operatorname{SL}_2(\mathbb{R})$ is a cocompact arithmetic subgroup. They show that $\norm{\phi}_\infty \ll \lambda^{1/4 - 1/24 + \varepsilon}$ for an $L^2$-normalised Hecke-Maaß form $\phi$.

Another parameter that one could choose is the \emph{volume} of $X$. This is particularly interesting in the arithmetic case and is reminiscent of the level aspect in the subconvexity problem of $L$-functions (indeed, the two problems are very much related in methodology and numerology). One thoroughly studied example is the family of non-compact spaces $X_0(N) := \Gamma_0(N) \backslash \mathfrak{h}^2$ of volume $N^{1 + o(1)}$, where $\Gamma_0(N)$ is the Hecke congruence subgroup of level $N$. 
For an example of an application, the corresponding `convexity' bound on average (for squarefree $N$) was used in \cite{AU} to compare the Arakelov and the Poincaré metrics on $X_0(N)$. A great amount of work was dedicated to achieving sub-baseline bounds in more and more general settings,  for example in \cite{BH}, \cite{templier}, \cite{HT}, \cite{saha-hybrid}, \cite{edgar}. See the introduction of \cite{saha-hu} for a more complete set of references with the corresponding bounds.

In general, the level aspect seems to factorise into the case of squarefree level and the case of powerful level, in particular prime powers. The latter is called the \emph{depth} aspect and is amenable to techniques from $p$-adic analysis that are not available in the squarefree case. To describe an example, let $N_1$ denote the smallest positive integer such that $N \mid N_1^2$. Note that $N_1 = N$ if $N$ is squarefree, yet $N_1 \asymp \sqrt{N}$ if $N$ is a high prime power. As an example of a subconvex bound in the case of $X_0(N)$ and the interplay between squarefree and powerful levels, it was shown in \cite{saha-hybrid} that
\begin{displaymath}
\norm{\phi}_\infty \ll N^{-2/6 +\varepsilon} N_1^{1/6} \lambda^{5/24 + \varepsilon},
\end{displaymath}
for $\phi$ a newform. This is also an example of a \emph{hybrid} bound. Such a bound is uniform in \emph{both} the spectral and the volume aspect.

Closer to the topic of these notes is the case of cocompact surfaces, where the arithmetic subgroup is given by the norm $1$ units $\mathcal{O}^1$ of an order $\mathcal{O}$ in an indefinite division quaternion algebra $A$ over $\mathbb{Q}$. The volume of $\mathcal{O}^1 \backslash \mathfrak{h}^2$ is approximately equal to the squareroot of the discriminant of $\mathcal{O}$ (see Section \ref{sec-volume}). This discriminant is made up of the discriminant of the algebra and the so-called level of $\mathcal{O}$. More precisely, we choose a fixed maximal order $\mathcal{O}_m$ containing $\mathcal{O}$ and define $N = [\mathcal{O}_m: \mathcal{O}]$, which we call the level, and note that $\disc(\mathcal{O}) = N^2 \cdot \disc(\mathcal{O}_m) = N^2 \cdot \disc(A)$. Most results in the literature only give bounds in the parameter $N$, but the discriminant of $A$, and therefore the covolume of $\mathcal{O}$, can get arbitrarily large independently of the level. A hybrid bound that is uniform in the full volume aspect was proved by Blomer and Michel \cite{blomer_michel} for the related case of totally definite quaternion algebras.

For Eichler orders of indefinite quaternion algebras, the ($p$-adic) local bound $\norm{\phi}_\infty \ll N^{-1/2+ \varepsilon} N_1^{1/2}$ was shown by Marshall \cite{marshall} for $\phi$ a newform. For squarefree $N$ this corresponds to the bound $N^\varepsilon$, which was first improved by Templier \cite{templier}, who obtained the bound $N^{-1/24 + \varepsilon}$ for general $N$. Saha \cite[Theorem 3]{saha} combines these two bounds to $N^{-11/24 + \varepsilon} N_1^{5/12}$ for newforms. Saha remarks in Section 1.6 that the argument should also provide a non-trivial hybrid bound but the details do not seem to be in print.

The depth aspect was recently improved in \cite{saha-hu} for newforms. If the level is a prime power $p^n$, then Hu and Saha obtain the bound $p^{n(5/24 + \varepsilon)}$. 

It should be noted at this point that the bounds above that include the number $N_1$ are proven using newform theory and, in particular, assuming that the order is Eichler, as stated. Saha and Templier also prove results for general orders (e.g.\ \cite[Theorem 1]{saha}), thus not assuming that the forms are newforms, but only that they are Hecke eigenfunctions. Templier \cite[Theorem 2]{templier} explicitly assumes $\phi$ is a newform, but the argument does not use newform theory; it does, however, use the structure of an Eichler order. As is apparent in the main results of this paper, described in the next section, our focus is rather on proving bounds in the case of general orders and we therefore never use newforms or the number $N_1$.

The sup-norm problem has also been pushed in the last decade to the case of higher-rank groups, such as $\operatorname{GL}(n)$ for $n > 2$. For example, Blomer and Maga \cite{BMn} prove that if $\phi$ is a Hecke-Maaß cusp form for the Hecke congruence group $\Gamma_0(N) \leq \operatorname{SL}_n(\mathbb{Z})$, then
\begin{displaymath}
\norm{\phi|_\Omega}_\infty \ll_\Omega N^\varepsilon \lambda^{\frac{n(n-1)}{8}(1 - \delta_n)},
\end{displaymath} 
for some fixed compact set $\Omega \subset \mathfrak{h}^n = \operatorname{SL}_n(\mathbb{R}) / \operatorname{SO}(n)$ and effectively computable $\delta_n > 0$. The local bound \eqref{eq-local} in this case is $n(n-1)/8$.

In higher-rank, the only other bound related to the volume aspect (that is available to the author) is a bound in the depth aspect given by Hu \cite{hu}. The result is stated only for automorphic forms corresponding to minimal vectors, which seem to be more suitable for the $p$-adic analysis of the depth aspect. 

\subsection{The main results and methods}

The purpose of this article is to establish the first non-trivial hybrid sup-norm bounds for Hecke-Maaß forms in \emph{unbounded} rank. Note that the uniformity in the volume (not just the level, but also the discriminant of the algebra) is a new feature even in the degree 2 case.

Before stating the main theorem, recall that a \emph{locally norm-maximal} order is an order $\mathcal{O}$ such that $\nr(\mathcal{O}^\times_p) = \mathbb{Z}_p^\times$ for all primes $p$, where $\mathcal{O}_p$ is the $p$-adic completion of $\mathcal{O}$. This a technical assumption in the theorems below and we refer to Remark \ref{rem-generalisation} for a way to remove it.

\begin{theorem}\label{thm-main}
Let $p \geq 3$ be a prime and $A$ a central division algebra of degree $p$ over $\mathbb{Q}$ that is split over $\mathbb{R}$. Let $\mathcal{O} \subset A$ be a locally norm-maximal order of covolume $V := \operatorname{vol} \mathcal{O}^1 \backslash \mathfrak{h}^p$. If $\phi$ is an $L^2$-normalised Hecke-Maaß form on $\mathcal{O}^1 \backslash \mathfrak{h}^p$ with large eigenvalue $\lambda$, then
\begin{equation}\label{eq-bound-p}
\norm{\phi}_\infty \ll \lambda^{\frac{p(p-1)}{8} - \delta_1 + \varepsilon} V^{- \delta_2 +\varepsilon},
\end{equation}
where the savings can be taken to be $\delta_1 = (16p^3)^{-1}$ and $\delta_2 = (8p^3(p-1))^{-1}$,
and the implied constant depends on $p$ and $\varepsilon$.
\end{theorem}

The so-called convexity bound is given by setting $\delta_1 = \delta_2 = 0$ in the exponents of \eqref{eq-bound-p}, so that we obtain subconvex bounds in both aspects simultaneously. The savings in the exponents also have the advantage of being explicit and polynomial in the degree. In fact, the proof shows a marginally stronger bound in the spectral aspect, but we have simplified the exponent for aesthetic reasons.

There are certain assumptions that can be removed in the statement of Theorem \ref{thm-main}, as well as in the theorems below. These include allowing automorphic forms that \emph{transform under characters} and orders that are \emph{not locally norm-maximal}, as noted in Remark \ref{rem-generalisation} after discussing the methods of proof.

The reason we only work over $\mathbb{Q}$ in Theorem \ref{thm-main} is explained in Remark \ref{rem-fields}. It is no loss of generality, but in fact the only field relevant to our problem for the space $\mathfrak{h}^p$ with $p > 2$. For the case $p=2$, we prove a bound over totally real number fields.

\begin{theorem}\label{thm-main-nf}
Let $A$ be an indefinite division quaternion algebra over a totally real number field $F$. Let $\mathcal{O} \subset A$ be a locally norm-maximal order of covolume $V := \operatorname{vol} \mathcal{O}^1 \backslash \mathfrak{h}^p$. If $\phi$ is an $L^2$-normalised Hecke-Maaß form on $\mathcal{O}^\times \backslash \mathfrak{h}^p$ with large eigenvalue $\lambda$, then
\begin{equation}
\norm{\phi}_\infty \ll \lambda^{\frac14 - \delta_1 + \varepsilon} \cdot V^{-\delta_2 + \varepsilon},
\end{equation}
where the savings can be taken to be $\delta_1 = 1/120$ and $\delta_2 = 1/30$, and the implied constant depends on $F$ and $\varepsilon$.
\end{theorem}

Considering the literature on the indefinite quaternion algebra case, the novelty of the bound is the uniformity in the discriminant of the algebra, and thus in the full volume of the quotient, while previous bounds only considered the level of the order. It should also be noted that this is the first hybrid bound appearing in print (accessible to the author). However, as a compromise for being hybrid, we should remark that our bounds in Theorem \ref{thm-main-nf} are weaker than those of \cite{templier} (i.e.\ $N^{-1/24}$, where we recall that $V \asymp_A N^{1 + o(1)}$) or \cite{IS} (i.e. $\lambda^{1/4-1/24}$) when isolating the relevant aspect.

In the case of quaternion algebras over the rational numbers, a sup-norm bound uniform in the full volume of the hyperbolic surface, in particular uniform in the discriminant, was also recently and independently achieved by Khayutin, Nelson, and Steiner \cite{KNS}. Among many impressive results, their theorem in the corresponding setting of Theorem \ref{thm-main-nf} is that $\norm{\phi}_\infty \ll_{\lambda, \varepsilon} V^{-1/4 + \varepsilon}$. This is a major improvement over other bounds in the literature, however only for Eichler orders and newforms, and without uniformity in the spectral aspect. Their proof uses different novel methods which, though very powerful in degree 2, do not seem to generalise easily to the higher rank situation.

To give the proof of our main theorems some historical context, we note that it shares some strategies with Section 6 of \cite{templier} and Section 2.4 of \cite{saha}, the latter being inspired by the former. In turn, the proof in \cite{templier} was inspired by the work of Silberman and Venkatesh \cite{SV} on quantum unique ergodicity, which has many similarities to the sup-norm problem. As \cite{SV} treats more generally division algebras of prime degree, it seems only natural that the previous sup-norm problem arguments should extend to this setting, and this is achieved in these notes. The main new idea for obtaining hybrid bounds is a systematic use of conjugation invariant functions, such as the characteristic polynomial, as we explain below.

We observe that going beyond prime degrees seems to require a new general strategy, as was also noted by Silberman and Venkatesh. The present paper provides a first step in this direction and the methods suffice for treating certain orders of Eichler type (see the discussion in Section \ref{sec-composite}), in arbitrary odd degree. We work with orders $\mathcal{O}$ of type $\mathcal{O}_0(N)$, by which we mean that, at unramified primes $p$, the completion $\mathcal{O}_p$ is of the form
\begin{displaymath}
	\mathcal{O}_0(N)_p = \left\{ \gamma \in M_n(\mathbb{Z}_p) \mid \text{last row of } \gamma \equiv (0, \ldots, 0, \ast)  \operatorname{mod} N \mathbb{Z}_p \right\},
\end{displaymath}
up to conjugation.

\begin{theorem}\label{thm-eichler}
Let $n \geq 3$ be an odd integer and $A$ a central division algebra of degree $n$ over $\mathbb{Q}$ that is split over $\mathbb{R}$. Let $\mathcal{O} \subset A$ be an order of type $\mathcal{O}_0(N)$ and let $V := \operatorname{vol} \mathcal{O}^1 \backslash \mathfrak{h}^n$ be its covolume. If $\phi$ is an $L^2$-normalised Hecke-Maaß form on $\mathcal{O}^1 \backslash \mathfrak{h}^n$ with large eigenvalue $\lambda$, then
\begin{equation}
\norm{\phi}_\infty \ll_A \lambda^{\frac{n(n-1)}{8} - \delta_1 + \varepsilon} V^{- \delta_2+\varepsilon},
\end{equation}
where the savings can be taken to be $\delta_1 = (8n^3)^{-1}$ and $\delta_2 = (4n^3(n-1))^{-1}$, and the implied constant depends on $n$, $\varepsilon$, and the discriminant of $A$.
\end{theorem}

Comparing the bound above with the one in Theorem \ref{thm-main}, the reader may remark that the savings are stronger by a factor of two. This is because, for the special orders to which it applies, the argument can handle counting in the spectral and level aspect at the same time, whilst for Theorem \ref{thm-main} we interpolate two different bounds, whence the halving of the savings.

Generalising Theorem \ref{thm-eichler} to handle any type of orders only requires improving the counting argument in the volume aspect. The uniform counting argument in the spectral aspect, given in Section \ref{sec-spec-count}, is valid for any orders, at least in odd degree.

The following discussion of the methods applies generally to all three theorems stated above. As usual in the treatment of the sup-norm problem, the argument starts with an amplified pretrace formula. We embed the form $\phi$ into a basis of Hecke-Maaß forms $(\phi_j)$ for $L^2(\mathcal{O}^1 \backslash \mathfrak{h}^p)$ and we spectrally expand an automorphic kernel with respect to this basis. This leads to an equality between a weighted sum (the spectral side) of the form
\begin{displaymath}
\sum_j A_j |\phi_j(z)|^2
\end{displaymath}
and a sum over elements in sets constructed from the group $\mathcal{O}^1$ (the geometric side).

To choose an amplifier essentially means to find suitable non-negative weights $A_j$ so that the contribution of $\phi$ is large and that of the other forms is little, hopefully negligible. These are constructed using Hecke eigenvalues.
This is a problem in analysis and combinatorics (or real and $p$-adic analysis), and was solved for example in \cite{BM4} for the groups $\operatorname{PGL}_n(\mathbb{R})$. Restricting to unramified places, we are also able to use the amplifier of Blomer and Maga.

After choosing an amplifier, we can then drop all but one terms and obtain a bound for $|\phi_j(z)|^2$ in terms of a sum over certain elements in $\mathcal{O}$, determined by the amplifier, which turns into a counting problem. We count elements $\gamma \in \mathcal{O}$ of norms up to some parameter $L$, such that the distance between $z$ and $\gamma z$ is small. This is usually done quite explicitly in the non-compact case of congruence subgroups of $\operatorname{SL}_2(\mathbb{Z})$. In our case, we rely on the very rigid structure of division algebras of prime degree.

More precisely, we may assume that the elements we are counting lie in a \emph{proper} subalgebra of $A$ at the cost of upper bounds for the parameter $L$. Here we make crucial use of the degree being prime (and nowhere else in an essential way). In this case, the only proper subalgebras of $A$ are fields, where we have better techniques available. In particular, it suffices to count ideals and units with certain conditions in the ring of integers, and the resulting number of elements is essentially best-possible. Thus, the bound in the theorem is dictated by how large we can take $L$ to be. We note that the scarceness of subalgebras in the prime degree case is also the reason why Silberman and Venkatesh prove their results in this setting (see Section 1.3 in \cite{SV}).

There are several difficulties to be overcome in the more general setting of this paper. One of them is that the counting problem we obtain when bounding $|\phi_j(z)|^2$ depends on the point $z$. This point only needs to vary along a compact region, but to make the bounds uniform in the volume, we need to have a good grasp on how the counting depends on $z$. In this paper we use conjugation invariant functions to completely remove $z$ from the counting problem.

Another difficulty is finding a way to incorporate the discriminant of the algebra, which is not visible in previous approaches. The crucial step is at Lemma \ref{lemma-proper}, where we deduce that a certain subspace spanned by elements in our order is proper inside the whole algebra. This is done using the determinant method, that is, by proving lower and upper bounds for a well-chosen determinant which contradict each other if the determinant is non-zero and certain parameters are small enough. The new input is the use of a special type of matrix whose determinant is divisible by the discriminant. This is a particularly useful choice of matrix, since it's entries are given by certain traces, which are conjugation invariant. This is a first example of the strategy mentioned in the previous paragraph.

The other conjugation invariant function we extensively make use of is the norm. It is essential for the uniform counting in the spectral aspect and for obtaining bounds in composite degrees. Our arguments not only obtain the properness of the algebra generated by the elements we are counting, but we prove its commutativity directly in certain cases. This is based on two simple but powerful observations: that commutators $k_1 k_2 - k_2 k_1$ of special orthogonal matrices in arbitrary odd degree are singular; and that commutators of elements in orders of $\mathcal{O}_0(N)$-type (e.g. matrices in $\Gamma_0(N)$) have determinant divisible by $N$. These observations, together with the fact that the only element in a division algebra with norm $0$ is $0$, solve the counting problem in the spectral and level aspect simultaneously for $\mathcal{O}_0(N)$-type orders and in the spectral aspect for any order. We refer to Section \ref{sec-spec-count} and Section \ref{sec-composite} for more details.

Another problem worth mentioning concerns the counting argument for units in a field (as explained above, we must also count ideals, but this is easily done by using unique factorisation). In this article, we handle this by bounding the possibilities for the characteristic polynomial of these elements, which again is conjugation invariant. Since we are counting in a (commutative) field, this automatically bounds the number of units. This argument is valid in any degree and can be employed whenever, by other arguments, one can restrict counting to a commutative algebra.


Finally, we remark that the counting argument makes use of the discriminant $\operatorname{disc}(\mathcal{O})$ in lieu of the volume of $\mathcal{O}^1 \backslash \mathfrak{h}^p$. It is well-known that these two parameters are essentially equal, as already mentioned in this introduction. This is indeed readily available for quaternion algebras, for which there is extensive literature available. Nevertheless, the author was not able to find a direct reference for the required formulae, at least for the more general setting in this article. We generalise the quaternion algebra argument to show that $\operatorname{vol}(\mathcal{O}^1 \backslash \mathfrak{h}^p) = \operatorname{disc}(\mathcal{O})^{1/2+o(1)}$ in section \ref{sec-volume}. The calculations may implicitly be present in other works on zeta functions of division algebras, but since this fact might be useful in particular in the theory of automorphic forms, we provide the details here to serve as reference.

\begin{remark}\label{rem-generalisation}
	The statements of the theorems presented in this introduction can be generalised by working adelically. Although this paper considers automorphic forms as classical objects, as in the theory of quantum chaos, both approaches are common in the literature (e.g.\ \cite{BM4} is written classically and \cite{saha} adelically).

	For example, taking $\phi$ in the theorems stated above to be an adelic form, one could drop the assumption that $\mathcal{O}$ is locally norm maximal (see also Remark \ref{rem-hecke}) and one could allow $\phi$ to generate a finite dimensional representation under the action of the corresponding adelisation of the unit group $\mathcal{O}^1$. Indeed, this is the approach taken in \cite{saha}, Theorem 1. Following the argument there, one would perform amplification using only unramified primes and the counting problem would be the same as in the classical case. Some details regarding the classical and adelic formulations of the theorems, in particular for automorphic forms transforming under a finite character of $\mathcal{O}^1$, are given in Section 1.3 of \cite{saha}. 

	The counting argument being the main novelty in this article, the classical formulation was chosen here to reduce technicalities and size. Some details on adelisation are given in Section \ref{sec-amp} to aid with understanding the classical Hecke algebra.
\end{remark}

\subsection*{Notation}
We recall the Vinogradov notation $f(x) \ll g(x)$ for two functions $f, g$, meaning that $f(x) \leq C \cdot g(x)$, at least for large enough $x$, for some $C > 0$ called the implied or implicit constant. We shall often use the notation $g = h + O(\delta)$ for matrices $g, h$, meaning that each coefficient of $g - h$ is $O(\delta)$. 
An expression of the form $f(x) = x^{o(1)}$ is to be intepreted as $f(x) \ll_\varepsilon x^\varepsilon$ and $x^{-\varepsilon} \ll_\varepsilon f(x)$ for any $\varepsilon > 0$, where the implied constant in both bounds can depend on $\varepsilon$.
Also, we sometimes work with more general degrees $n \in \mathbb{N}$ for the division algebra $A$, and restrict where necessary to prime degrees $p$.

\subsection*{Acknowledgements}
I would like to thank Valentin Blomer for introducing me to this topic and for his constant support. For many enlightening discussions on these topics, I thank my colleagues Edgar Assing and Bart Michels. For their kind help with understanding their work, I would also like to thank Abhishek Saha and John Voight. Finally, the readability and accuracy of the paper were greatly improved by the anonymous referee's comments and suggestions, for which I am grateful.

\section{Division algebras and arithmetic subgroups of $\operatorname{SL}(n, \mathbb{R})$}

Let $A$ be a central division algebra of degree $n$ over $\mathbb{Q}$.  Let $\mathcal{O} \subset A$ be an order, i.e.\ a subring with $1$ that is a full $\mathbb{Z}$-lattice. Suppose that $A$ splits over $\mathbb{R}$, meaning that there is an embedding $A \hookrightarrow M_n(\mathbb{R})$. For an element $x \in A$, the reduced norm $\operatorname{nr} (x)$ and the reduced trace $\operatorname{tr} (x)$ are given by the determinant and the trace, respectively, of its image under this embedding. The group $\mathcal{O}^1 = \{ \gamma \in \mathcal{O} :  \operatorname{nr}(\gamma) = 1\}$ now embeds into $\operatorname{SL}_n(\mathbb{R})$ as a cocompact arithmetic lattice (see \cite{morris}, Proposition 6.8.9). In particular, denoting the symmetric space of $\operatorname{SL}_n(\mathbb{R})$ by
\begin{displaymath}
\mathfrak{h}^n = \operatorname{SL}_n(\mathbb{R}) / \operatorname{SO}(n),
\end{displaymath}
then the quotient $\mathcal{O}^1 \backslash \mathfrak{h}^n$ is compact.

Note that for $n$ odd, $A$ splits automatically over $\mathbb{R}$. Indeed, by the Albert-Brauer-Hasse-Noether theorem, all central simple algebras of finite degree over a number field are cyclic, meaning that they contain a strictly maximal subfield that is a Galois extension of $\mathbb{Q}$ of degree $n$ (for background on these statements and the following, see \cite{pierce}, Theorem 18.6 and sections 13.1 through 13.3). The strictly maximal subfield of $A$ splits $A$ and is Galois of odd degree over $\mathbb{Q}$, so it must be totally real, in particular contained in $\mathbb{R}$.

In the special case $n = 2$, $A$ is called a quaternion algebra and we may replace the ground field $\mathbb{Q}$ by any totally real number field $F$. Let $[F:\mathbb{Q}] = n$ and denote by $\mathfrak{o}_F$ the ring of integers of $F$. For $A$ to be split over $\mathbb{R}$, we assume that there is an embedding $\sigma_0 \in \hom(F, \mathbb{R})$ such that \(A \otimes_{\sigma_0} \mathbb{R} \cong M_2(\mathbb{R})\). For all other embeddings $\sigma_0 \neq \sigma \in \hom(F, \mathbb{R})$ assume that \(A \otimes_{\sigma} \mathbb{R} \cong \mathcal{H}(\mathbb{R})\), where $\mathcal{H}(\mathbb{R})$ is the Hamilton quaternion algebra. We may view $A$ as embedded (diagonally) into $A_\infty \cong M_2(\mathbb{R}) \times \mathcal{H}(\mathbb{R})^{n-1}$, and similarly for the norm $1$ elements,
\begin{displaymath}
A^1 \hookrightarrow \operatorname{SL}_2(\mathbb{R}) \times \operatorname{SO}(3)^{n-1}.
\end{displaymath}
We use $\phi_0$ to denote the projection onto the first component $M_2(\mathbb{R})$ and $\phi_i$, $i = 1, \ldots, n-1$, to denote the projections onto the Hamiltonian components. 

Generalising our setting, let $\mathcal{O}$ be an $\mathfrak{o}_F$-order. By restriction of scalars, the projection $\phi_0(\mathcal{O}^1) \subset \operatorname{SL}_2(\mathbb{R})$ of the group of units of reduced norm 1 onto the split component gives a cocompact arithmetic lattice.

\begin{remark}\label{rem-fields}
For $\phi_0(\mathcal{O}^1)$ to be a cocompact arithmetic subgroup in the split component $\operatorname{SL}_2(\mathbb{R})$, it is important that the other components, in this case all isomorphic to $\operatorname{SO}(3)$, are compact (see the definition of an arithmetic group in \cite{morris}, Definition 5.1.19). If $A$ is a central division algebra over a number field $F \neq \mathbb{Q}$ and $\operatorname{deg}(A) = n > 2$, this is not possible any more. 

Indeed, the process of restriction of scalars requires us to embed $A^1$ into the product of its completions at all infinite places. Now a central simple algebra over $\mathbb{R}$ is isomorphic to a matrix algebra over a division algebra by Wedderburn's theorem. By a theorem of Frobenius (see \cite{pierce}, Corollary 13.1 c) these are either matrix algebras over $\mathbb{R}$ or over the Hamiltonians $\mathcal{H}$. Since $n > 2$, the group of norm 1 units in these algebras cannot be compact any more. Thus, the number field case in higher degree gives rise to non-compact lattices, which are outside the scope of this article.
\end{remark}

It will be useful later to note that the tower rule holds for division
algebras (also called skew fields). More precisely, the notion of vector space
over a division algebra and its dimension is the same as for commutative
fields. If $A' \subset A$ is a subalgebra, then
$A$ may be viewed as a vector space over $A'$, where $A'$ acts by
multiplication from the left (or from the right, according to taste). We
denote the dimension of the vector space by $\dim_{A'} A$ as usual.

Let now $A_1 \subset A_2 \subset A_3$ be division algebras. Then 
\begin{equation}\label{eq-tower}
\dim_{A_1} A_3 = \dim_{A_2} A_3 \cdot \dim_{A_1}A_2
\end{equation}
holds and is proven as in the commutative
case  (see \cite[Proposition 3.1.1]{cohn}). Thus, if $A$ is a finite dimensional division algebra over
$\mathbb{Q}$, then the dimension over $\mathbb{Q}$ of any subalgebra of $A$
must divide $\dim_{\mathbb{Q}} A$. Moreover, if $A$ is central, then any
subfield of $A$ must have dimension over $\mathbb{Q}$ dividing the degree of
$A$ (see \cite{pierce}, Corollary 13.1 a).

\subsection{The volume approximation}\label{sec-volume}

For simplicity, we first assume that the ground field is $\mathbb{Q}$ and quote the relevant results for quaternion algebras over number fields at the end of the section.

Let $\mathcal{O}_m$ be a maximal order in $A$ containing $\mathcal{O}$. Because of their lattice structure, it is useful to work with the index $[\mathcal{O}_m : \mathcal{O}] $, which we call the \emph{level} of $\mathcal{O}$ in $\mathcal{O}_m$. Yet the volume of $\mathcal{O}^1 \backslash \mathfrak{h}^n$, the relevant parameter in our sup-norm problem, is given by the volume of $\mathcal{O}_m^1 \backslash \mathfrak{h}^n$ and the multiplicative index $[\mathcal{O}_m^1 : \mathcal{O}^1] $. Fortunately the two indices are related in an explicit way. For our purposes (and because the exact formulae would involve too many cases in general), it suffices to prove that they are approximately equal. The proper equalities obtained in the proof can be used together with the machinery of zeta functions and Tamagawa numbers to produce a formula for the volume of $\mathcal{O}^1 \backslash \mathfrak{h}^n$ (as in \cite{voight}, 39.2.8).

\begin{lemma}\label{lemma-index}
Let $A$ be a central simple algebra over $\mathbb{Q}$, but not a definite quaternion algebra. Let $\mathcal{O} \subset \mathcal{O}_m$ be two orders in $A$. If $[\mathcal{O}_m : \mathcal{O}]  = N$, then $[\mathcal{O}_m^1 : \mathcal{O}^1]  = N^{1+o(1)}$.
\end{lemma}

This lemma is a generalisation of the well-know fact that the index of the Hecke congruence group $\Gamma_0(N)$ in $\operatorname{SL}_2(\mathbb{Z})$ is 
\begin{displaymath}
N \prod_{p \mid N} (1 + 1/p).
\end{displaymath}
In this case, $A$ is the matrix algebra $M_2(\mathbb{Q})$, the maximal order is $M_2(\mathbb{Z})$ and $\mathcal{O}$ is the suborder of level $N$ of integral matrices with lower left entry divisible by $N$.

Before proving Lemma \ref{lemma-index}, we use it together with a formula for the covolume of a maximal order to approximate the covolume of $\mathcal{O}$ by the discriminant. Recall that the \emph{discriminant} of an arbitrary $\mathbb{Z}$-order $\mathcal{O}$ is defined as the ideal $\operatorname{disc}(\mathcal{O}) \subset \mathbb{Z}$ generated by the set
\begin{displaymath}
\{ \det  (\operatorname{tr}(x_i \cdot x_j))_{1 \leq i, j \leq n} \mid x_i \in \mathcal{O} \}.
\end{displaymath}
By abuse of notation, we also denote a positive generator of this ideal by $\operatorname{disc}(\mathcal{O})$. For more details on discriminants, see Section 10 in \cite{reiner}.

\begin{proposition}\label{prop-vol}
Let $A$ be a central simple algebra over $\mathbb{Q}$, but not a definite quaternion algebra, and let $\mathcal{O}$ be an order in $A$. Then $\operatorname{vol}(\mathcal{O}^1 \backslash \mathfrak{h}^n) = \operatorname{disc}(\mathcal{O})^{1/2+o(1)}$.
\end{proposition}
\begin{proof}
Let $\mathcal{O}_m$ be a maximal order containing $\mathcal{O}$. By Theorem 3.7 in \cite{kleinert}, we can approximate
\begin{displaymath}
\vol(\mathcal{O}^1_m \backslash \mathfrak{h}^n) = \disc(\mathcal{O}_m)^{\frac{1}{2} + o(1)}.
\end{displaymath}
By Lemma 15.2.15 in \cite{voight}, we also have
\begin{displaymath}
\disc(\mathcal{O}) = [\mathcal{O}_m: \mathcal{O}]^2 \cdot \disc(\mathcal{O}_m).  
\end{displaymath}
Together with Lemma \ref{lemma-index}, we now obtain the claimed approximation since $\vol(\mathcal{O}^1 \backslash \mathfrak{h}^n) = \vol(\mathcal{O}_m^1 \backslash \mathfrak{h}^n) \cdot [\mathcal{O}_m^1 : \mathcal{O}^1]$.
\end{proof}

The proof of Lemma \ref{lemma-index} generalises the argument for quaternion algebras in \cite{voight}, Lemma 26.6.7, which in turn follows an argument of Körner. We provide here full details for the sake of completeness. 

The first ingredient is the strong approximation theorem (see Kneser's article in \cite{boulder}), which allows us to reduce the statement to a local one. We denote by
$A_p = A \otimes \mathbb{Q}_p$ and $\mathcal{O}_p =\mathcal{O} \otimes
\mathbb{Z}_p$ the completions at a prime $p$. For all but finitely many primes
$p$, the completion $A_p$ is split, i.e.\ $A_p \cong \mathcal{M}_n
(\mathbb{Q}_p)$ (see Proposition 18.5 coupled with Corollary 17.10.a in {\cite{pierce}}).
Additionally, for all but finitely many primes $p$, the completion
$\mathcal{O}_p$ is a maximal order of $A_p$ (see Lemma 10.4.4 in
{\cite{voight}}). In particular, at these primes we have $\mathcal{O}_{m,p} = \mathcal{O}_p$. The primes where equality does not hold will be referred to as \emph{ramified}. 

We embed $\mathcal{O}$ diagonally into $\hat{\mathcal{O}} = \prod_p \mathcal{O}_p$ and, similarly, $ \mathcal{O}^1 $ into $\hat{\mathcal{O}}^1 = \prod_p \mathcal{O}_p^1$, where $p$ runs over all prime numbers. Then strong approximation implies that $\mathcal{O}^1$ is dense in $\hat{\mathcal{O}}^1$ (see \cite{voight}, Corollary 18.5.14, and more generally \cite{kleinert}, Theorem 4.4). Explicitly, if $S$ a set of finite places, $a_p \in \mathcal{O}^1_{p}$ and $t_p$ for each $p \in S$, then we can find $x \in \mathcal{O}^1$ such that
\begin{equation*}
x \equiv a_v  \pmod{ p^{t_v} \mathcal{O}_p} \quad  (p \in S).
 \end{equation*}

\begin{lemma}\label{lemma-local-ind}
	For two orders $\mathcal{O} \subset \mathcal{O}_m$ as above, the level and the index of the unit groups can be computed locally, that is,
	\begin{equation*}
		[\mathcal{O}_m : \mathcal{O}] = \prod_p [\mathcal{O}_{m,p} : \mathcal{O}_p] \quad \text{and} \quad [\mathcal{O}_m^1 : \mathcal{O}^1] = \prod_p [\mathcal{O}_{m,p}^1 : \mathcal{O}_p^1].
	\end{equation*}
\end{lemma}

\begin{proof}
	Note first that the products contain only finitely many factors not equal to $1$, as in the remarks above. Next, we start the proof for the unit groups. The claim follows by showing that the map
	\begin{equation*}
		\mathcal{O}_m^1 / \mathcal{O}^1 = \prod_p \mathcal{O}_{m,p}^1 / \mathcal{O}_p^1
	\end{equation*}
	is bijective.
	
	Injectivity follows by noting that $\bigcap_p \mathcal{O}_p = \mathcal{O}$. Surjectivity follows by strong approximation. Indeed, let $(a_p) \in \prod_p \mathcal{O}_{m,p}^1$. Choose an integer $N$ such that $N \mathcal{O}_{m,p} \subset p \mathcal{O}_p$ for all ramified primes $p$. Strong approximation supplies us with an element $b \in \mathcal{O}_m^1$ such that $b \in a_p + N\mathcal{O}_{m,p}$. Thus $b = a_p \cdot u_p$, where $u_p \in 1 +  N\mathcal{O}_{m,p}$, so that $u_p \in \mathcal{O}_p^1$.
	
	The proof for the factorisation of the level is similar, where the corresponding strong approximation theorem is the Chinese Remainder Theorem.
\end{proof}

In the following we work with the localised orders at a prime $p$, which we suppress in notation for simplicity. We now remove the condition on the norm to work with the full group of units. We have the short exact sequence
\begin{equation*}
	0 \rightarrow \mathcal{O}^1 \rightarrow \mathcal{O}^\times \rightarrow \operatorname{nr}(\mathcal{O}^\times) \rightarrow 0,
\end{equation*}
and similarly for $\mathcal{O}_m$. By defining a non-canonical bijection\footnote{Note that the groups in question are not abelian and not necessarily normal, so that we cannot apply the snake lemma directly.} $\mathcal{O}_m^1 / \mathcal{O}^1 \times \operatorname{nr}(\mathcal{O}_m^\times) / \operatorname{nr}(\mathcal{O}^\times) \rightarrow \mathcal{O}_m^\times / \mathcal{O}^\times$, we obtain that
\begin{equation*}
	| \mathcal{O}_m^1 / \mathcal{O}^1 | \cdot | \operatorname{nr}(\mathcal{O}_m^\times) / \operatorname{nr}(\mathcal{O}^\times) | = | \mathcal{O}_m^\times / \mathcal{O}^\times |.
\end{equation*}

\begin{lemma}
For $\mathbb{Z}_p$-orders $\mathcal{O} \subset \mathcal{O}_m$, we have
\begin{displaymath}
[ \mathcal{O}_m^\times : \mathcal{O}^\times ] = [ \mathcal{O}_m : \mathcal{O} ] \cdot p^{o(1)}.
\end{displaymath}
\end{lemma}

\begin{proof}
The proof starts as in Lemma 26.6.7 in \cite{voight}.
Let $n$ be such that $p^n \mathcal{O}_m \subset p \mathcal{O}$. Note that $1 + p\mathcal{O} \subset \mathcal{O}^\times$ by the convergence of the geometric series. We now have
\begin{displaymath}
[ \mathcal{O}_m^\times : \mathcal{O}^\times ] = \frac{ [\mathcal{O}_m^\times : 1 + p\mathcal{O}_m ] \cdot [1 + p \mathcal{O}_m: 1 + p^n \mathcal{O}_m] } { [\mathcal{O}^\times : 1 + p\mathcal{O}] \cdot [1 + p \mathcal{O} : 1 + p^n \mathcal{O}_m]}.
\end{displaymath}
For $\alpha, \beta \in 1 + p \mathcal{O}$, we have $\alpha \beta^{-1} \in 1 + p^n \mathcal{O}_m$ if and only if $\alpha - \beta \in p^n \mathcal{O}_m$, so that
\begin{displaymath}
[1 + p \mathcal{O} : 1 + p^n \mathcal{O}_m] = [p \mathcal{O}: p^n \mathcal{O}_m] = [\mathcal{O}: p^{n-1}\mathcal{O}_m],
\end{displaymath}
and similarly for $\mathcal{O}_m$. If follows that
\begin{displaymath}
\frac{[1 + p \mathcal{O}_m: 1 + p^n \mathcal{O}_m]} {[1 + p \mathcal{O} : 1 + p^n \mathcal{O}_m]} = \frac{[\mathcal{O}_m : p^{n-1}\mathcal{O}_m]} {[\mathcal{O}: p^{n-1}\mathcal{O}_m]} = [\mathcal{O}_m : \mathcal{O}].
\end{displaymath}

To further compute the factors $[\mathcal{O}^\times : 1 + p \mathcal{O}]$, we employ the strategy in \cite{voight}, Lemma 24.3.12, of introducing the Jacobson radical $\operatorname{rad}\mathcal{O} =: J$. We have $p \mathcal{O} \subset J$ and there is an integer $r$ such that $J^r \subset p \mathcal{O} $ (see \cite{reiner}, Theorem 6.13), which we assume to be minimal. Thus $1 + J \subset \mathcal{O}^\times$ and we obtain a filtration
\begin{displaymath}
\mathcal{O}^\times \supset 1 + J \supset 1 + J^2 \supset \ldots \supset 1 + p \mathcal{O} \supset 1 + J^r.
\end{displaymath}
Being kernels, all subgroups are normal inside their parent groups. It follows that 
\begin{displaymath}
[\mathcal{O}^\times : 1 + p \mathcal{O}] =  |\mathcal{O}^\times / 1 + J | \cdot |1 + J / 1 + J^2 | \cdots |1 + J^{r-1} / 1 + p \mathcal{O}|.
\end{displaymath}

On the additive side, we also have a filtration $\mathcal{O} \supset J \supset \ldots \supset p \mathcal{O}$ and the quotients $\mathcal{O}/J, J/J^2, \ldots, J^{r-1}/p \mathcal{O}$ are $\mathbb{F}_p$-algebras. If $R$ is the rank of $\mathcal{O}$, then
\begin{displaymath}
p^R = [\mathcal{O} : p \mathcal{O}] = [\mathcal{O}: J] [J : J^2] \cdots [J^{r-1} : p \mathcal{O}].
\end{displaymath}

We now reduce the multiplicative indices to the additive ones. Indeed $1 + J / 1 + J^2 \cong J / J^2$, and similarly for all powers of $J$, and $1 + J^{r-1} / 1 + p \mathcal{O} \cong J^{r-1}/p \mathcal{O}$, since $J^{2(r-1)} \subset p \mathcal{O}$ (at least for $r > 1$; the case $r = 1$ is simpler and can be done directly). Therefore,
\begin{displaymath}
|1 + J / 1 + J^2 | = |J / J^2|, \ldots, |1 + J^{r-1} / 1 + p \mathcal{O}| = |J^{r-1} / p \mathcal{O}|.
\end{displaymath}
Thus,
\begin{displaymath}
[\mathcal{O}^\times : 1 + p \mathcal{O}] = p^R \frac{ [\mathcal{O}^\times / 1 + J] } { [\mathcal{O} : J] }.
\end{displaymath}

Now one can easily see that $\mathcal{O}^\times / 1 + J \cong (\mathcal{O} / J)^\times$. The reason for working with the Jacobson radical is that $\mathcal{O} / J$ is a semisimple $\mathbb{F}_p$-algebra, meaning that 
\begin{displaymath}
\mathcal{O} / J \cong M_{d_1}(A_1) \times \cdots \times M_{d_l}(A_l),
\end{displaymath}
for some finite division algebras $A_i$ over $\mathbb{F}_p$. Since finite division algebras are fields by Wedderburn's theorem, one can check by counting that $|GL_{d_i}(A_i)| = |M_{d_i}(A_i)|^{1- o(1)} $ (this can be made precise, but the approximation is sufficient for our purposes). 

Since $\mathcal{O}$ and $\mathcal{O}_m$ have the same rank, it follows that
\begin{displaymath}
\frac{ [\mathcal{O}_m^\times : 1 + p \mathcal{O}_m] } { [\mathcal{O}^\times : 1 + p \mathcal{O}] } = p^{o(1)}.
\end{displaymath}
This finishes the proof.
\end{proof}

For the groups of norms, we note that $\mathbb{Z}_p^\times \subset \mathcal{O}^\times$, and so $(\mathbb{Z}_p^\times)^n \leq \operatorname{nr}(\mathcal{O}^\times) \leq \mathbb{Z}_p^\times$. Now $[ \mathbb{Z}_p^\times : \mathbb{Z}_p^\times)^n] \ll n$ by Korollar 5.8 in \cite{neukirch}. This contributes to the global index by $n^{\omega(N)} \ll_n d(N) \ll N^\varepsilon$, where $\omega(N)$ is the number of different primes dividing the level $N$, $d(N)$ is the number of divisors of $N$, and $\varepsilon$ is any positive real number. This completes the proof of Lemma \ref{lemma-index}.

A similar argument, taking into account all of the different embeddings, would also apply for division algebras over number fields. Since we are only interested in quaternion algebras for totally real number fields, we simply note an approximate version of Main Theorem 39.1.8 in \cite{voight} (it is given for locally norm-maximal orders; see Remark 39.1.12 for the general case, which is the same observation we make in the previous paragraph). Let $F$ be a totally real number field with ring of integers $\mathfrak{o}_F$. In this case, we recall that the discriminant $\disc(\mathcal{O})$ is an $\mathfrak{o}_F$-ideal and we denote by $N_F(\disc(\mathcal{O}))$ its norm.

\begin{proposition}\label{prop-vol-nf}
Let $A$ be an indefinite quaternion algebra over $F$ and let $\mathcal{O}$ be an order in $A$. Then $\vol (\mathcal{O}^1 \backslash \mathfrak{h}^2) = N_F (\disc(\mathcal{O}))^{1/2+o(1)}$, where the implicit constant depends on $F$.
\end{proposition}

\section{The amplified pretrace formula}\label{sec-amp}

In this section we again assume that the ground field is $\mathbb{Q}$. The additional technicalities involved in the case of quaternion algebras over number fields are described in Section \ref{sec-quat}.

The space of automorphic forms $L^2 ( \mathcal{O}^1 \backslash \mathfrak{h}^n ) $ has a discrete decomposition, admitting a basis of Hecke-Maaß forms $(\phi_j)_{j \in \mathbb{N}}$, that is, eigenfunctions of the algebra of invariant differential operators and of the Hecke algebra (described below). Denote the spectral parameters of each form $\phi_j$ by $\mu_j \in \mathfrak{a}^\ast_\mathbb{C}$, where $\mathfrak{a}^\ast_\mathbb{C}$ is the complexified dual of the Lie algebra of the diagonal torus $A$ in $\operatorname{SL}_n(\mathbb{R})$. Recalling the main goal of this paper we note that for bounding an individual automorphic form $\phi$ with spectral parameter $\mu$, we may assume that $\phi$ is part of the basis $(\phi_j)$. 

For setting up the pretrace formula, we follow the notation in \cite{BMn}. In particular, to any $\lambda \in \mathfrak{a}^\ast_\mathbb{C}$ we attach the quantity $D(\lambda)$ defined in \cite[(1.2)]{BMn}.  To ease notational clutter due to the fact that the discriminant of $\mathcal{O}$ is denoted by $D$, we put $S(\lambda) := D(\lambda)$ in this paper. As in Section 2, ibid., we denote $\mu^\ast := \Re \mu \in \mathfrak{a}^\ast$ and assume that $\norm{\mu^\ast}$ is sufficiently large.  If $\lambda_\phi$ is the Laplace eigenvalue of $\phi$, then 
\begin{equation}\label{eq-eigenvalue}
S(\mu^\ast) \ll 1 + \lambda_\phi^{n(n-1)/4},
\end{equation}
as in \cite[(2.3)]{BMn}. The bounds we obtain in the remainder of the paper use $S(\mu^\ast)$ instead of the eigenvalue since they are slightly improved this way (we only stated the main theorems using the eigenvalue for simplicity), but also to simplify some exponents.

As usual we denote $K = \operatorname{SO}(n)$. For a function $f \in C_c^\infty(K \backslash G / K)$, the pretrace formula states that
\begin{displaymath}
\sum_{j \in \mathbb{N}} \tilde{f}(\mu_j) \phi_j(z) \overline{\phi_j(z')} = \sum_{\gamma \in \mathcal{O}^1} f(z^{-1} \gamma z'),
\end{displaymath}
for $z, z' \in G$, where $\tilde{f}$ is the spherical transform of $f$. 

To calibrate the pretrace formula for our distinguished function $\phi$ with spectral parameter $\mu$, Blomer and Maga (e.g. see \cite{BMn}, Section 3)  show that we can find $f_\mu \in C_c^\infty (K \backslash G / K)$ such that the spherical transform $\tilde{f}$ satisfies $\tilde{f}_\mu (\lambda) \geq 0$ for all possible spectral parameters $\lambda$ and $\tilde{f}_\mu (\mu) \geq 1$. 

\begin{remark}
We may use the same test function $f_\mu$ as Blomer and Maga, since the relevant spectral parameters of automorphic forms for the division algebra $A$ depend only on the archimedean completion, more precisely on the Lie group $\SL_n(\mathbb{R})$. These are included in the set defined in (2.2) of \cite{BMn} by the general theory of joint eigenfunctions and spherical functions on symmetric spaces (see \cite{helgason}, in particular Chapter IV, Section 1, 2, and 8). Generally, the choice of test function $f_\mu$ is a local problem at the archimedean place, while applying Hecke operators as we do below is a choice of test function at the finite places (the reader comfortable with the adelic theory of automorphic forms could see Section 4.1 of \cite{saha} for the adelic treatment of amplification).
\end{remark}

In fact, using bounds of Blomer and Pohl \cite[Sect. 6]{BP}, we can also assume certain decay properties of $f_\mu$. More precisely, we can assume the diameter of the support of $f_\mu$ to be bounded by any positive constant depending on the (fixed) degree $n$ as we prefer. Moreover, if $d$ denotes the invariant distance function on $\mathfrak{h}^n$, then we have the bound 
\begin{equation} \label{eq-spec-bound}
f_\mu(g) \ll S(\mu^\ast) (1 + S(\mu^\ast)^{\frac{2}{n(n-1)}} \cdot d(g, 1))^{-1/2},
\end{equation}
which is easily implied by \cite[(2.4)]{BMn}.

To further amplify the contribution of $\phi$ in the pretrace formula, Blomer and Maga also construct a general amplifier using Hecke operators (see \cite{BM4}, Section 6) for $\operatorname{SL}_n(\mathbb{Z})$. This amplifier applies in our situation as well, as long as we only use unramified places. To explain this statement, we sketch below some facts about the Hecke algebra, for which we note our assumption that the ground field is $\mathbb{Q}$.

First, we define the group $U_{\mathcal{O}}$ as
\begin{displaymath}
	U_{\mathcal{O}} = \text{GL}_n^+ (\mathbb{R}) \times \prod_p \mathcal{O}_p^{\times}
\end{displaymath}
Note that 
\begin{displaymath}
	\mathcal{O}^1 = U_{\mathcal{O}} \cap A^\times.
\end{displaymath}

Next, as can be seen from Lemma \ref{lemma-local-ind} for instance, the primes $p$ that divide $D := \disc(\mathcal{O})$ are exactly the primes at which $A$ ramifies or $\mathcal{O}_p$ is not maximal, and we call these primes \emph{ramified}. 
%
%
%
%

From now on we assume that $\mathcal{O}$ is \emph{locally norm-maximal}, meaning that $\nr(\mathcal{O}^\times_p) = \mathbb{Z}_p^\times$ for all primes $p$. This assumption, explained in Remark \ref{rem-hecke}, is satisfied, for example, by any maximal orders or intersection of two maximal orders, i.e.\ Eichler orders. We define the semigroup $S_{\mathcal{O}}$ inside the adelisation $A_{\mathbb{A}}^{\times}$ by
\begin{displaymath}
	S_{\mathcal{O}} = \left( \text{GL}_n^+ (\mathbb{R}) \times \prod_p S_p \right) \cap A_{\mathbb{A}}^{\times},
\end{displaymath}
where $S_p = \{ \alpha \in \mathcal{O}_p : \operatorname{nr} (\alpha) \neq 0 \} $ for $p \nmid D$ and $S_p = \mathcal{O}_p^\times $ for $p \mid D $. This distinction means that we only consider the \emph{unramified} Hecke algebra. Finally, let
\begin{displaymath}
	\Delta_{\mathcal{O}} = S_{\mathcal{O}} \cap \mathcal{O}.
\end{displaymath}

We can now define the \emph{classical} Hecke algebra $R(\mathcal{O}^1, \Delta_\mathcal{O})$, which is generated by double cosets of the form $\mathcal{O}^1 \xi \mathcal{O}^1$, where $\xi \in \Delta_{\mathcal{O}}$, and similarly the \emph{adelic} Hecke algebra $R(U_\mathcal{O}, S_\mathcal{O})$. For more details, see \cite{miyake}, Sections 2.7 and 5.3. 

The adelic point of view is advantageous since we automatically obtain a factorisation of $R(U_\mathcal{O}, S_\mathcal{O})$ as the tensor product $\bigotimes_p R(\mathcal{O}_p^\times, S_p)$ of the local Hecke algebras. Fortunately in our case, there is essentially nothing lost in translation between the classical and the adelic Hecke algebra (both unramified). Indeed, they are isomorphic under the simple correspondence $\mathcal{O}^1 \xi \mathcal{O}^1 \mapsto U_\mathcal{O} \xi U_\mathcal{O}$. This can be seen by carefully applying the argument in the proof of Theorem 5.3.5 in \cite{miyake}. The proof makes crucial use of approximation theorems.

\begin{remark}\label{rem-hecke}
The property of being locally norm-maximal implies that the idelic quotient defined by $\mathcal{O}$ has only one connected component. In particular, the dictionary between classical automorphic forms and adelic forms is simpler. For many works in the literature (see for instance the use of Eichler orders in \cite{templier}, and \cite{SV}, Remark 6.3.1), this is a common ``cosmetic'' assumption on orders.

Indeed, our counting arguments in prime degree, the key new ideas in this paper, simply make no use of this assumption. Removing it is merely a matter of working directly with adelic automorphic forms and the full unramified Hecke algebra. This is the approach in \cite{saha}, where the generalisation of Templier's argument for Eichler orders \cite{templier} to arbitrary orders is transparent. 
	
In fact, even working classically, the pretrace inequality \eqref{eq-pt-ineq} below and Remark \ref{rem-L} are still true for any order $\mathcal{O}$, at least conditionally on a suitable Riemann hypothesis. We sketch this technicality in this paragraph, which should be read preferably after going through the proof of Theorem \ref{thm-main}. More precisely, to adjust the definition of the Hecke algebra, one may need to impose additionally that $S_p = \mathcal{O}_p^\times $ for all primes $p$ that are not $n$-th powers modulo $D$ (supposing for simplicity that $(D, n) = 1$). This is the ``worst-case'' scenario, since we have the inclusions $(\mathbb{Z}^\times_p)^n \subset \nr(\mathcal{O}_p^\times) \subset \mathbb{Z}^\times_p$. Under this assumption, the proof of  Theorem 5.3.5 in \cite{miyake} goes through. We are now left with a potentially smaller set of primes $\mathcal{P}$, defined just before \eqref{eq-pt-ineq}, where we also assume that these primes are $n$-th powers modulo $D$. An application of the generalised Riemann hypothesis and the Chebotarev density theorem would then imply that $|\mathcal{P}| \gg_n L^{1-\varepsilon} D^{-\varepsilon}$, even for $L$ a small power of $D$, as taken in the proof of the main theorem.
\end{remark}

Now at unramified primes $p$, the local Hecke algebra $R(\mathcal{O}_p^\times, S_p)$ is isomorphic to the Hecke algebra of $\operatorname{GL}_n(\mathbb{Z}_p)$. Therefore, we can use the same Hecke operators as Blomer and Maga.

For $m \in \mathbb{Z}$, let
\begin{displaymath}
\mathcal{O}(m) := \{ \gamma \in \mathcal{O} \mid \operatorname{nr}(\gamma) = m \}.
\end{displaymath}
Let $L > 5$ be a parameter and $\mathcal{P}$ be the set of primes in $[L, 2L]$ that are unramified. We have the pretrace inequality (see \cite{BMn}, (2.5))
\begin{equation}\label{eq-pt-ineq}
|\mathcal{P}|^2 \cdot |\phi(z)|^2 \ll |\mathcal{P}|\cdot S(\mu^\ast) + \sum_{\nu = 1}^n \sum_{l_1, l_2 \in \mathcal{P}} \frac{1}{L^{(n-1)\nu}} \sum_{\gamma \in \mathcal{O}(l_1^\nu l_2^{(n-1)\nu})} |f_\mu(z^{-1} \tilde{\gamma} z) |,
\end{equation}
where $\tilde{\gamma} = \gamma / \operatorname{nr}(\gamma)^{1/n} \in \operatorname{SL}_n(\mathbb{R})$. Note that the determinantal divisors in \cite{BM4} do not seem to easily translate into our setting, yet the norm does, as one can easily check using the explicit isomorphism between the classical and adelic Hecke algebras above. Although these additional conditions are very important in the work of Blomer and Maga, we are able to solve the counting problem described below using only the condition on the norm.

Since $f_\mu$ has compact support, let $0< \rho \ll 1$ be such that $f_\mu(g) = 0$ if $d(g, 1) > \rho$, where $d$ is the invariant distance function on $\mathfrak{h}^n$. Using the bound $f_\mu \ll S(\mu^\ast)$ from \eqref{eq-spec-bound}, we may obtain an explicit bound for $\phi(z)$ from the pretrace inequality by counting the number of elements $\gamma \in \mathcal{O}(m)$ such that $d(z, \tilde{\gamma} z) < \rho$. We correspondingly define in general
\begin{displaymath} 
\mathcal{O} (m ; z, \delta) = \{ \gamma \in \mathcal{O}: \operatorname{nr}( \gamma ) = m, d (z, \tilde{\gamma} z) < \delta \}.
\end{displaymath} 

\begin{remark}\label{rem-supp}
We note that the compact support of $f_\mu$ can be assumed to be small enough in terms of the degree $n$, i.e. $\rho \ll_n 1$ with an implicit constant of our choice, since we are allowing all implicit constants to depend on $n$. This follows from a quick inspection of the technique in \cite{BP}.
\end{remark}

A more careful use of \eqref{eq-spec-bound} gives us a saving in the spectral aspect, i.e. in $S(\mu^\ast)$, at least if $d(z, \tilde{\gamma} z) > \delta$ for $\delta > 0$ large enough in terms of $S(\mu^\ast)$. For all other $\gamma$ appearing in the sum, we trivially bound $f_\mu(z^{-1} \tilde{\gamma} z) \ll S(\mu^\ast)$ as above. We therefore split the sum using the parameter $\delta$ as above and obtain the inequality
\begin{multline}\label{eq-pt-ineq-split}
|\mathcal{P}|^2 \cdot |\phi(z)|^2 \ll S(\mu^\ast) \bigl(
|\mathcal{P}| + \sum_{\nu = 1}^n \sum_{l_1, l_2 \in \mathcal{P}} \frac{1}{L^{(n-1)\nu}} \#\mathcal{O}(l_1^\nu l_2^{(n-1)\nu}; z, \delta) \\
+ S(\mu^\ast)^{\frac{-1}{n(n-1)}} \cdot \delta^{-\frac{1}{2}} \sum_{\nu = 1}^n \sum_{l_1, l_2 \in \mathcal{P}} \frac{1}{L^{(n-1)\nu}} \#\mathcal{O}(l_1^\nu l_2^{(n-1)\nu}; z, \rho)  \bigr).
\end{multline}

\begin{remark}\label{rem-L}
Note that $|\mathcal{P}| \gg L^{1 - \varepsilon} \cdot D^{-\varepsilon}$, at least for $L$ large enough. This follows from the prime number theorem and because the number of ramified primes we leave out is bounded by $\tau(D) \ll D^\varepsilon$. 
\end{remark}

To obtain a hybrid bound, we need to count elements in $\mathcal{O} (m ; z, \delta)$ uniformly in $\disc (\mathcal{O})$ and $S(\mu^\ast)$. To deduce anything about the sup-norm of $\phi$, the counting must also be done uniformly in $z$, at least in a fundamental domain for $\mathcal{O}^1$. Though compact, this fundamental domain grows with the volume. 

In \cite{saha}, the problem of uniformity in $z$ was resolved by having a counting argument that only depends on the index of $\mathcal{O}$ inside a maximal order $\mathcal{O}_m$. Saha was then able to conjugate $z$ into a fixed fundamental domain for $\mathcal{O}_m$ and work with a conjugated order with the same index. The implicit constants in the bounds would then possibly depend on the particular choice of fundamental domain for the maximal order.

Our argument does not require a choice of maximal order and the bounds are uniform on the whole generalised upper half plane. This is done by counting using only traces and norms, which are conjugation invariant. 

More precisely, the condition $d(z, \tilde{\gamma} z) \ll \rho$ implies that $z^{-1} \tilde{\gamma} z = k + O(\rho)$ for some $k \in \operatorname{SO}(n)$, at least for $\rho \ll 1$ small. Indeed, using the Cartan decomposition, we can write any $g \in \SL_n(\mathbb{R})$ as $g = k_1 \exp(C(g)) k_2$, where $k_1, k_2 \in \operatorname{SO}(n)$ and $C(g)$ is diagonal with vanishing trace. Then $d(g,1) = \norm{C(g)}_2$, where we view $C(g)$ as a vector in $\mathbb{R}^n$. The claim follows by writing the exponential as a power series.

On the new condition, as an example, applying the trace directly already provides a bound on $\tr(\tilde{\gamma})$ where any dependence on $z$ is completely removed, noting that orthogonal matrices are bounded. This conjugation invariant approach is used throughout the counting argument. Another example is given in the subsection below, where we derive the so-called convexity bound for the sup-norm automorphic forms in certain cases.


%
%
%

\subsection{The baseline bound}
The benchmark bound for the sup-norm of automorphic forms that one seeks to improve can be usually obtained by using the pretrace formula without the amplifier. It follows easily from the properties of $f_\mu$ that
\begin{displaymath}
|\phi(z)|^2 \leq \sum_{\gamma \in \mathcal{O}^1} f_\mu(z^{-1} \gamma z) \ll S(\mu^\ast) \cdot \# \{ \gamma \in \mathcal{O}^1 : z^{-1} \gamma z = k + O(\rho) \text{, for some } k \in \operatorname{SO}(n) \}.
\end{displaymath}

Suppose that $n$ is odd. Then any degree $n$ special orthogonal matrix must have $1$ as an eigenvalue, meaning that $\det (k - 1) = 0$. Thus, if $\gamma \in \mathcal{O}^1$ and $z^{-1} \gamma z = k + O(\rho)$, then we may subtract the identity matrix and apply the determinant to obtain
\begin{displaymath}
\nr(\gamma -1) = \det( k - 1 + O(\rho)) = O_n(\rho).
\end{displaymath}
As mentioned in Remark \ref{rem-supp}, we can take $\rho$ as small as we wish in terms of $n$. Since $\gamma - 1 \in \mathcal{O}$, it follows by integrality of the norm that $\nr(\gamma - 1) = 0$. Since $A$ is a division algebra, this implies that $\gamma = 1$. Therefore, there is only one term appearing on the geometric side of the pretrace formula formula and using the inequality \eqref{eq-spec-bound} we obtain
\begin{displaymath}
|\phi(z)| \ll S(\mu^\ast)^{\frac12},
\end{displaymath}
where the implied constant depends at most on $n$. This is the convexity bound (recall also the bound \ref{eq-eigenvalue}).

This strategy cannot succeed when $n$ is even, if only for the simple observation that $-1$ lies in $\mathcal{O}^1$. Some ad-hoc arguments involving the classification of motions suffice for the case $n = 2$, but the author was not able to find a general argument for all even $n$. This is part of the reason why this paper mainly deals with algebras of odd degree, besides quaternion algebras.

\section{Counting in the discriminant aspect}\label{sect-disc-count}


For simplicity, we first describe the counting argument over $\mathbb{Q}$ and adjust it where necessary for quaternion algebras over number fields in the next subsection. From now on, assume that $\operatorname{deg}(A) = p \geq 3$, a prime. Let $\delta$ be a positive real number, which we assume to be uniformly bounded, e.g. by $1$. As explained in the previous section, we are interested in bounding the cardinality of
\begin{displaymath} 
\mathcal{O} (m ; z, \delta) = \{ \gamma \in \mathcal{O}: \operatorname{nr}( \gamma ) = m, d (z, \tilde{\gamma} z) < \delta \},
\end{displaymath}
where $\tilde{\gamma} = \gamma / \operatorname{nr}(\gamma)^{1/p}$. The condition $d (z, \tilde{\gamma} z) < \delta$ is equivalent to $z^{-1} \tilde{\gamma} z = k + O(\delta)$ for some $k \in \operatorname{SO}(p)$.

To motivate the following lemmata, we recall the tower rule \eqref{eq-tower} for division algebras. Especially for prime degree $p$, this severely restricts the possible dimensions of subalgebras in $A$. If one can show that the subalgebra generated by the elements we are counting is proper, then the tower rule drastically reduces the dimension of the counting problem, automatically. In our case, the subalgebra will actually be commutative, which is crucial in our argument. To show properness in the first place, we use a version of the determinant method, for which we need good control  over a basis of a vector space.

\begin{lemma} \label{lemma-alg-in-vs}
The $\mathbb{Q}$-algebra generated by $\bigcup_{1 \leq m \leq L} \mathcal{O} (m ; z, \delta)$ is contained in the $\mathbb{Q}$-vector space spanned by $\bigcup_{1 \leq m \leq L^{2p-2}} \mathcal{O} (m ; z, (2p-2) \delta)$.
\end{lemma}

\begin{proof}
By the tower rule, a subalgebra of $A$ is of the form $\mathbb{Q}$, $\mathbb{Q} (x)$, or $\mathbb{Q} (x, y)$, where $x, y \in A$ and $\mathbb{Q}(x,y)$ is the smallest algebra containing both $x$ and $y$. Indeed, if a subalgebra contains a non-rational element $x$, then it contains $\mathbb{Q}(x)$, which must have dimension $p$ over $\mathbb{Q}$ (the characteristic polynomial has degree $p$). If it contains another element not in $\mathbb{Q}(x)$, say $y$, then it contains $\mathbb{Q}(x,y)$. The tower rule implies now that $A = \mathbb{Q}(x,y)$. The algebra $\mathbb{Q} (x, y)$ is generated as a vector space by monomials of degree at most $2p-2$.

Now if $a_j \in \bigcup_{1 \leq m \leq L} \mathcal{O} (m ; z, \delta)$ for $j = 1, \ldots, 2p-2$, then the reduced norm of $\prod a_j$ is at most $L^{2p-2}$ and, by the triangle inequality, $d (z, \prod \tilde{a_j} \cdot z) < (2p-2) \delta$. The order structure ensures that $\prod a_j$ lies in $\mathcal{O}$.
\end{proof}

For the next lemma denote $\disc (\mathcal{O}) := D$.

\begin{lemma} \label{lemma-proper}
The $\mathbb{Q}$-vector space spanned by $\bigcup_{1 \leq m \leq 	L^{2p-2}} \mathcal{O} (m ; z, (2p-2) \delta)$ is proper, i.e.\ not equal to $A$, if $L \ll D^{1/4p(p-1) - \varepsilon}$, where the implicit constant depends only on $p$ and $\delta$.
\end{lemma}
\begin{proof}
Let $\gamma_1, \ldots, \gamma_{p^2}$ be elements of $\bigcup_{1 \leq m \leq L^{2p-2}} \mathcal{O} (m ; z, (2p-2) \delta)$. In this case $\operatorname{nr}(\gamma_i \gamma_j) \ll L^{4(p-1)}$ and we have that $d(\gamma_i \gamma_j z, z) \leq d(\gamma_i z, z) + d(\gamma_j z, z) < 4(p-1)\delta$, by the triangle inequality. In particular 
\begin{displaymath}
\operatorname{tr}(\gamma_i \gamma_j) \ll_{\delta, p} L^{4(p-1)/p},
\end{displaymath}
by applying the trace to the equation 
$z^{-1} \gamma_i \gamma_j z = \nr(\gamma_i \gamma_j)^{1/p} ( k + O_p(\delta) )$,
with some $k \in \operatorname{SO}(p)$.

Consider now $s = \det(\operatorname{tr}(\gamma_i \gamma_j)_{i,j})$. Recall that $D$ is the generator of the ideal in $\mathbb{Z}$ generated by $\{ \det\operatorname{tr}(x_i x_j) \mid x_i \in \mathcal{O}, i = 1,\ldots p^2\}$. Since $\gamma_i \in \mathcal{O}$ for all $i$, it follows that $D \mid s$. 

On the other hand, by using the bound above, we deduce that $s \ll L^{4p(p-1)}$. Thus if $L \ll D^{1/4p(p-1) - \varepsilon}$, then $s = 0$. By the non-degeneracy of the bilinear form given by the reduced trace, it follows that $\gamma_1, \ldots, \gamma_{p^2}$ are \emph{not} linearly independent.
\end{proof}

Thus, if $L$ is small enough, we can assume that we are counting matrices in a proper subalgebra of $A$, which must be $\mathbb{Q}$ or a field extension $E /\mathbb{Q}$ of degree $p$. This is where the use of $p$ being a prime is crucial. 

We are now counting certain elements in $\mathcal{O}_E$, the ring of integers of $E$, which certainly includes $\mathcal{O} \cap E$. We do so by counting ideals and units. Since the units in $\mathbb{Z}$ are only $\pm 1$, we can concentrate on the non-trivial extensions, which must have an infinite group of units, at least if $p>2$. It is important to note that the reduced norm and the reduced trace in $A$ of an element in a subfield $E \subset A$ such that $\dim_{\mathbb{Q}} E = p$ are the same as the number field norm, resp. trace of $E / \mathbb{Q}$ (see \cite[Sect. 16.2]{pierce}).

\begin{lemma} \label{lemma-units-field}
Let $E /\mathbb{Q}$ be a cyclic extension of degree $p$ that is a subfield of $A$ and let $\mathcal{O}$ be an order of $A$. The number of units $\xi \in \mathcal{O}^{\times} \cap E$ such that $d (z, \tilde{\xi} z) \leq \delta$ for a given $z \in \mathfrak{h}^p$ is $\ll p^p (1+\delta)^{p-1}$.
\end{lemma}

\begin{proof}
Let $\xi \in \mathcal{O}^{\times} \cap E$. Then $\xi \in \mathcal{O}_E^{\times}$, where $\mathcal{O}_E$ is the ring of integers of $L$, by integrality over $\mathbb{Z}$. Thus, $\operatorname{nr} (\xi) = N_{E /\mathbb{Q}} (\xi) = \pm 1$.

Next, the condition $d (z, \tilde{\xi} z) \leq \delta$ is equivalent to $\xi \in z B (\delta) z^{- 1}$, where $B (\delta)$ is a union of $\delta$-balls around all elements of $\text{SO} (p)$. Applying the trace, we find that $\operatorname{tr}(\xi) \ll p (1+ \delta)$. Since $\operatorname{tr} (\xi) \in \mathbb{Z}$ by integrality, we see that there are $\ll p (1+ \delta)$ possibilities for the value of $\operatorname{tr} (\xi)$.

We may apply the same reasoning to $\xi^j$ and derive that there are $\ll p (1 + j \delta)$ possibilities for the value of $\operatorname{tr} (\xi^j)$. Indeed, 
\begin{displaymath}
d (z, \tilde{\xi}^j z) \leq d (z, \tilde{\xi} z) + d (\xi z, \tilde{\xi}^{j} z) = d (z, \tilde{\xi} z) + d (z, \tilde{\xi}^{j-1} z) \leq j \delta,
\end{displaymath}
inductively. Note also that $\xi^j \in \mathcal{O}^{\times} \cap E$ since $\mathcal{O}$ is closed under multiplication.

Now the characteristic polynomial of $\xi$ is
\begin{displaymath}
X^p - \operatorname{tr} (\xi) X^{p-1} + \frac{1}{2}\left[ \operatorname{tr} (\xi)^2 + \operatorname{tr} (\xi^2)\right] X^{p-2} + \ldots \pm \det (\xi).
\end{displaymath}
By Newton's identities, each coefficient is determined by the values of $\operatorname{tr} (\xi^j)$ for certain $j$. By the bounds above, there are only $\ll \prod_{j=1, \ldots {p-1}} p(1+j \delta) \ll [p(1+\delta)]^{p-1}$ polynomials that are satisfied by a unit $\xi$ as in the statement of the lemma. Since each polynomial can have at most $p$ different roots, the proof is finished.
\end{proof}

\begin{lemma}\label{lemma-bound-in-field}
Let $E \subset A$ be a field of degree $p$ over $\mathbb{Q}$. Then for a any $z \in \mathfrak{h}^p$ and any positive integer $m$ we have
\begin{displaymath}
\mathcal{O}_E (m ; z, \delta) \ll_p \tau (m)^{p-1} \cdot (1 + \delta)^{p-1} .
\end{displaymath}
\end{lemma}

\begin{proof}
Let $\gamma \in \mathcal{O}_E$ with $\operatorname{nr} (\gamma) = N_{E /\mathbb{Q}} (\gamma) = m$. Up to units, there are only $\tau (m)^{p-1}$ elements of $\mathcal{O}_E$ with norm $m$. Indeed, a principal ideal is determined by its generator up to units and the norm of the ideal is equal to the norm of the generator. Since ideals factorise uniquely into prime factors, we only need to count prime ideals.

Above each rational prime, there are at most $p$ prime ideals of $\mathcal{O}_E$. Therefore, if $q^v \mid m$ for a prime $q$, then we need to choose at most $p$ numbers $a_1, \ldots, a_p \in \mathbb{Z}_{\geq 0}$ such that $a_1 + \ldots + a_p = v$ to determine an ideal of norm $q^v$. The number of such tuples is $v+p-1$ choose $p-1$, that is $<<v^{p-1}$. Thus, there are at most 
\begin{displaymath}
\ll \prod_{q^v \| m} v^{p-1} = \tau (m)^{p-1}
\end{displaymath}
ideals of norm $m$.

Now if $\gamma \in \mathcal{O}_E (m ; z, \delta)$ and $\xi \gamma \in \mathcal{O}_E (m ; z, \delta)$ for some unit $\xi \in
\mathcal{O}_E^{\times}$, then 
\begin{displaymath}
d (z, \tilde{\xi} z) \leq d (z, \tilde{\xi} \tilde{\gamma} z) + d (\tilde{\xi} \tilde{\gamma} z, \xi z) = d (z, \tilde{\xi} \tilde{\gamma} z) + u (z, \tilde{\gamma} z) \leq 2
\delta.
\end{displaymath}
Thus, we finish the proof by counting such units using Lemma \ref{lemma-units-field}.
\end{proof}

Putting everything together and recalling the divisor bound $\tau(m) \ll m^\varepsilon$ we obtain the following proposition.

\begin{proposition}\label{prop-disc-count}
Let $m \ll D^{\frac{1}{4p(p-1)} - \varepsilon}$ with implicit constant as in Lemma \ref{lemma-proper} and let $z \in \mathfrak{h}^p$. Then $\# \mathcal{O} (m ; z, \delta) \ll_p m^{\varepsilon}$, where the implicit constant depends only on $\varepsilon$, $\delta$, and $p$.
\end{proposition}

\section{Counting in the spectral aspect}\label{sec-spec-count}

Let $\mu$ be the spectral
parameter of $\phi$ and denote $S \asymp D (\mu^{\ast})$ as in
\cite[(1.2)]{BMn}.

\subsection{The counting argument for small $\delta$}

We are interested in bounding the cardinality of
\begin{displaymath} 
\mathcal{O} (m ; z, \delta) = \{ \gamma \in \mathcal{O}: \operatorname{nr}( \gamma ) = m, d (z, \tilde{\gamma} z) = O (\delta) \},
\end{displaymath}
where $\tilde{\gamma} = \gamma / \operatorname{nr}(\gamma)^{1/p}$. The condition $d (z, \tilde{\gamma} z) = O (\delta)$ is equivalent to $z^{-1} \tilde{\gamma} z = k + O(\delta)$ for some $k \in \operatorname{SO}(p)$.

As opposed to the discriminant aspect, we can now gain savings in the pretrace inequality by using $\delta$ as a parameter.  
For $p \geq 3$, instead of showing that the $\mathbb{Q}$-algebra generated by the elements we are counting is a proper algebra and then deducing commutativity, we directly show that it must be a commutative field, at least if $\delta$ is small enough.

For the following lemma, assume that $p \geq 3$ is any odd integer, not necessarily prime.
\begin{lemma} \label{lemma-comm-delta}
The $\mathbb{Q}$-algebra generated by $\bigcup_{1 \leq m \leq L} \mathcal{O} (m ; z, \delta)$ is commutative, i.e.\ a field, if $\delta \ll  L^{-2 - \varepsilon}$, where the implicit constant depends only on $p$.
\end{lemma}
\begin{proof}
Let $\gamma_1, \gamma_2 \in \bigcup_{1 \leq m \leq L} \mathcal{O} (m ; z, \delta)$. A few applications of the triangle inequality, recalling that $d(\tilde{\gamma} z, z) = d(\tilde{\gamma}^{-1} z, z)$,\footnote{For instance, $d(\alpha^{-1} \beta z, z) \leq d(\alpha^{-1} \beta z, \alpha^{-1} z) + d(\alpha^{-1} z, z) = d(\beta z, z) + d(\alpha z, z)$ by invariance of the distance function.} show that $d(\gamma_1^{-1} \gamma_2^{-1} \gamma_1 \gamma_2 z, z) \ll \delta$. 
This implies that 
\begin{equation} \label{eq-commutator}
z^{-1} \gamma_1^{-1} \gamma_2^{-1} \gamma_1 \gamma_2 z = k + O(\delta),
\end{equation}
since $\nr(\gamma_1^{-1} \gamma_2^{-1} \gamma_1 \gamma_2) = 1$.

Recall that we assume $p \geq 3$ and note that $\det(k - 1) = 0$ for all $k \in \operatorname{SO}(p)$, since $1$ is certainly an eigenvalue of orthogonal matrices in odd degree. Therefore, subtracting $1$ from \eqref{eq-commutator} and taking the determinant gives
\begin{align*}
\nr( \gamma_1^{-1} \gamma_2^{-1} \gamma_1 \gamma_2 - 1) &= \det( z^{-1} (\gamma_1^{-1} \gamma_2^{-1} \gamma_1 \gamma_2 - 1) z ) \\
&= \det( k -1 + O(\delta)) = O_p(\delta).
\end{align*}
Multiplying this last equation by $\nr(\gamma_2 \gamma_1)$ implies that
\begin{equation}\label{eq-comm-bound}
\nr(\gamma_1 \gamma_2 - \gamma_2 \gamma_1) = \nr(\gamma_2 \gamma_1) \cdot \nr( \gamma_1^{-1} \gamma_2^{-1} \gamma_1 \gamma_2 - 1) = O_p(\delta \nr(\gamma_1) \nr(\gamma_2)).
\end{equation}

The commutator of two elements of $\mathcal{O}$ is again in $\mathcal{O}$ and so $\nr(\gamma_1 \gamma_2 - \gamma_2 \gamma_1) \in \mathbb{Z}$. Thus if $\delta \ll L^{-2 - \varepsilon}$, then it follows that $\nr(\gamma_1 \gamma_2 - \gamma_2 \gamma_1) = 0$, which implies that $\gamma_1 \gamma_2 = \gamma_2 \gamma_1$ since $A$ is a division algebra.
\end{proof}

The argument above makes crucial use of the fact that special orthogonal matrices in odd degree necessarily have $1$ as eigenvalue, which is not the case any more in even degree. To produce an alternative argument for $p = 2$, we shall return to the strategy employed for the discriminant aspect. This is done in Section \ref{sec-count-quat}.

Given the lemma above, we may assume that we are counting in a field $\mathbb{Q} \subset E$. This is done exactly as in Lemma \ref{lemma-bound-in-field}.

\begin{proposition}\label{prop-spec-count}
Let $A$ have odd degree $p \geq 3$ . If $\delta \ll m^{-2 - \varepsilon}$, then $\# \mathcal{O} (m ; z, \delta) \ll m^\varepsilon$.
\end{proposition}

\subsection{The counting argument for large $\delta$}

If $\delta = \rho \gg 1$ is as large as the diameter of the support of $f_\mu$ (see Section \ref{sec-amp}), we still need to bound the cardinality of $\mathcal{O} (m ; z, \delta)$ by a reasonable power of $m$ but with no dependence on the discriminant $D$. Since the amplifier is constructed only from unramified primes, as explained in Section \ref{sec-amp}, we may assume that $m$ is coprime to the discriminant of $A$.

\begin{proposition}\label{prop-spec-triv-count}
We have the bound $\# \mathcal{O} (m ; z, \rho) \ll m^{p-1+ \varepsilon}$, where the implicit constant depends only on $p$.
\end{proposition}
\begin{proof}
Let $\gamma \in \mathcal{O} (m ; z, \rho)$ and consider the principal ideal $\gamma \mathcal{O}_m$ in a maximal order $\mathcal{O}_m$ containing $\mathcal{O}$. We have $\nr(\gamma) \mathbb{Z} = \nr(\gamma \mathcal{O}_m)$ and if $\gamma \mathcal{O}_m = \gamma' \mathcal{O}_m$, then $\gamma \xi = \gamma'$ for some unit $\xi \in \mathcal{O}_m^\times$. It therefore suffices to count ideals with norm $m$ and bound the number of units $\xi$ as above.

By the local-global dictionary for ideals, it suffices to count locally. Since $m$ is coprime to the discriminant of $A$, we only need to count ideals of norm $q^e$ in $M_p(\mathbb{Z}_q)$ for primes $q$ and positive integers $e$. This is done for instance as in Lemma 26.4.1 in \cite{voight} and is a well-known computation in the theory of zeta functions of algebras. We use the fact that all ideals are principal and employ the theory of elementary divisors to find a generator of the ideal of lower triangular form, where the diagonal is given by $(q^{a_1}, \ldots, q^{a_p})$ for non-negative integers $a_1, \ldots, a_p$ such that $a_1 + \ldots a_p = e$. All entries on the column below $q^{a_j}$ are uniquely defined as elements of $\mathbb{Z}/q^{a_j} \mathbb{Z}$. An easy counting argument thus shows that there are
\begin{displaymath}
\sum_{a_1 + \ldots a_p = e} q^{(p-1) a_1 + (p-2) a_2 + \ldots a_{p-1}} \ll q^{e (p-1) (1 + \varepsilon)}
\end{displaymath}
ideals of norm $q^e$ in $M_p(\mathbb{Z}_q)$. 

This shows that there are $\ll m^{p-1}$ ideals of norm $m$ in $\mathcal{O}_m$, which is implicitly also a bound for the number of principal ideals. Now let $\xi \in \mathcal{O}^\times_m$ be a unit as in the first paragraph of this proof. By the triangle inequality, we find that $d(\xi z, z) \ll \rho$. Since $\nr(\xi) = 1$, this implies that $z^{-1} \xi z = k + O(\rho)$, for some $k \in \operatorname{SO}(p)$. 

Suppose $p$ is odd. If $\xi_1, \xi_2$ are two such units, then the same reasoning as before \eqref{eq-comm-bound} provides the bound
\begin{displaymath}
\nr(\xi_1 \xi_2 - \xi_2 \xi_1) = O_p(\rho).
\end{displaymath}
As observed in Remark \ref{rem-supp}, we may assume that $\nr(\xi_1 \xi_2 - \xi_2 \xi_1) < 1$. By integrality, this implies that $\xi_1 \xi_2 = \xi_2 \xi_1$. 

Therefore, we may reduce the counting problem to counting units in the maximal order of a field by restricting to the commutative algebra generated by these units. As proved in Lemma \ref{lemma-bound-in-field}, there are $\ll_p 1$ units with the required properties. Putting all bounds together proves the statement.
\end{proof}

\section{Proof of Theorem \ref{thm-main}}\label{sec-proof}
We are now ready to insert the counting results into the pretrace inequality \eqref{eq-pt-ineq-split}. Let the degree $n$ be an odd number, at least 3.

\subsection{The spectral aspect}
We take the parameter $\delta$ as large as possible to still have control over the counting problem, but also gain a saving in the third term of \eqref{eq-pt-ineq-split}. Since the largest norm coming up in the pretrace inequality is $l_1^n l_2^{(n-1)n} \ll L^{n^2}$, we must take $\delta \ll L^{-2n^2 - \varepsilon}$ to be able to apply Proposition \ref{prop-spec-count}.

With such a choice of $\delta$, the second term in \eqref{eq-pt-ineq-split} is bounded as
\begin{equation}\label{eq-2nd-term}
\sum_{\nu = 1}^n \sum_{l_1, l_2 \in \mathcal{P}} \frac{1}{L^{(n-1)\nu}} \#\mathcal{O}(l_1^\nu l_2^{(n-1)\nu}; z, \delta) \ll_n |\mathcal{P}|^2 L^{-(n-1) + \varepsilon}.
\end{equation}

As already mentioned, to obtain a saving from the third term, we take $\delta$ as large as possible, i.e. $\delta \asymp L^{-2n^2 - \varepsilon}$. Together with Proposition \ref{prop-spec-triv-count}, the third term in \eqref{eq-pt-ineq-split} can be bounded by
\begin{align*}
S(\mu^\ast)^\frac{-1}{n(n-1)} \cdot \delta^\frac{-1}{2} \sum_{\nu = 1}^n \sum_{l_1, l_2 \in \mathcal{P}} \frac{1}{L^{(n-1)\nu}} \#\mathcal{O}(l_1^\nu l_2^{(n-1)\nu}; z, \rho)  &\ll 
S(\mu^\ast)^\frac{-1}{n(n-1)} L^{n^2 + \varepsilon} |\mathcal{P}|^2 \sum_{\nu = 1}^n \frac{L^{(n-1)^2 \nu}} {L^{(n-1)\nu}} \\
&\ll_n S(\mu^\ast)^\frac{-1}{n(n-1)} |\mathcal{P}|^2 L^{n^2 + n(n-1)(n-2) + \varepsilon}.
\end{align*}

We introduce the two bounds above into the pretrace inequality and, recalling that $|\mathcal{P}| \gg L^{1 - \varepsilon} \cdot D^{-\varepsilon}$, where $D$ is the discriminant, we obtain
\begin{displaymath}
|\phi(z)|^2 \ll S(\mu^\ast) \left( L^{-1 + \varepsilon} D^\varepsilon + L^{-(n-1)+\varepsilon} D^\varepsilon + S(\mu^\ast)^\frac{-1}{n(n-1)} L^{n^2 + n(n-1)(n-2) + \varepsilon} \right).
\end{displaymath}
It is clear that the first term dominates the second. After solving an easy optimisation problem for the first and third terms, we arrive at the bound 
\begin{equation}\label{eq-spec-final}
|\phi(z)|^2 \ll S(\mu^\ast)  \cdot S(\mu^\ast)^{\frac{-1}{n(n-1)(n^2 + n(n-1)(n-2) + 1)}+\varepsilon} \cdot D^\varepsilon \ll S(\mu^\ast) \cdot S(\mu^\ast)^{\frac{-1}{n^4(n-1)}+\varepsilon}  \cdot D^\varepsilon,
\end{equation}
where we weaken the first bound to the last one simply for aesthetic reasons.

\subsection{The discriminant aspect}
Here we assume that $n = p$ is a prime. We cannot gain any saving in the discriminant aspect by using the $\delta$ parameter. Therefore, we set $\delta = \rho$ and ignore the last term in $\eqref{eq-pt-ineq-split}$. In this case, we take $L$ as large as Proposition \ref{prop-disc-count} allows. Taking into account the largest norms appearing in \eqref{eq-pt-ineq-split}, we may choose $L \asymp D^{\frac{1}{4n^3(n-1)}}$. Then the second term in the pretrace inequality can be bounded as in \eqref{eq-2nd-term}, so that we have 
\begin{equation}\label{eq-disc-final}
|\phi(z)|^2 \ll S(\mu^\ast) (L^{-1 + \varepsilon} + L^{-(n-1)+\varepsilon}) \ll S(\mu^\ast) \cdot L^{-1 + \varepsilon} \ll S(\mu^\ast) D^{\frac{-1}{4n^3(n-1)}+\varepsilon}.
\end{equation}

\subsection{The hybrid bound}

We interpolate between the two bounds \eqref{eq-spec-final} and \eqref{eq-disc-final} simply by multiplying them, so that
\begin{displaymath}
|\phi(z)| \ll |\phi(z)|^{1/2} \cdot |\phi(z)|^{1/2} \ll S(\mu^\ast)^{\frac12} S(\mu^\ast)^{\frac{-1}{4n^4(n-1)}+\varepsilon} D^{\frac{-1}{16n^3(n-1)}+\varepsilon}.
\end{displaymath}
Recalling \eqref{eq-eigenvalue}, this proves Theorem \ref{thm-main}.

\section{Quaternion algebras over number fields}\label{sec-quat}

The full classification of cocompact arithmetic subgroups requires us to also consider quaternion algebras over number fields and Templier \cite{templier} treats the counting problem in this more general setting (though only in the level aspect and only for Eichler orders). This case is slightly more technical and we treat the problem by applying the same ideas as above more carefully. We first recall the theoretical background.

Let $F$ be a totally real number field of degree $n$. We denote by $\mathfrak{o}_F$ its ring of integers and by $N_F$ the number field norm of $F/\mathbb{Q}$. Let $A$ be a division quaternion algebra over $F$ and assume that there is an embedding $\sigma_0 \in \hom(F, \mathbb{R})$ such that \(A \otimes_{\sigma_0} \mathbb{R} \cong M_2(\mathbb{R})\). For all other embeddings $\sigma_0 \neq \sigma \in \hom(F, \mathbb{R})$ assume that \(A \otimes_{\sigma} \mathbb{R} \cong \mathcal{H}(\mathbb{R})\), where $\mathcal{H}(\mathbb{R})$ is the Hamilton quaternion algebra. 

Now let $\mathcal{O}$ be an $\mathfrak{o}_F$-order and let $\mathfrak{D} = \disc(\mathcal{O}) \subset \mathfrak{o}_F$ be its discriminant. By abuse of notation, we also denote a generator of the discriminant ideal by $\mathfrak{D}$. Let $$D := |N_F(\mathfrak{D})|.$$

We may view $A$ as embedded in $A_\infty \cong M_2(\mathbb{R}) \times \mathcal{H}(\mathbb{R})^{n-1}$. We use $\varphi_0$ to denote the projection onto the first component $M_2(\mathbb{R})$ and $\varphi_i$, $i = 1, \ldots, n-1$, to denote the projections onto the Hamiltonian components.

Note that $\sigma_i(\operatorname{tr} (\gamma)) = \operatorname{tr} (\varphi_i(\gamma))$ for $\varphi_i$ the projection onto the $\sigma_i$ component. The trace on the left hand side refers to the quaternion trace and on the right hand side it refers to the usual matrix trace. Similarly, $\sigma_i (\operatorname{nr} (\gamma)) = \operatorname{nr}(\varphi_i(\gamma))$.

\subsection{Remarks on the amplifier}
We consider automorphic forms on $\mathfrak{h}^2$ invariant under the arithmetic group $\mathcal{O}^\times$, viewed as a subgroup of $\operatorname{PGL}_2(\mathbb{R})$.
\footnote{This is a slightly different subgroup than in the rest of the paper, where we take norm $1$ units. It is a technical assumption due to the fact that the units $\mathfrak{o}_F^\times$ (or even the totally positive units, depending on the setup) might not all be squares. This implies that the quotient $\mathcal{O}^1 \backslash \operatorname{SL}_2(\mathbb{R})$ might be larger than $\mathbb{R}^\times \mathcal{O}^\times \backslash \operatorname{GL}_2(\mathbb{R})$. The latter is the more natural one from the point of view of automorphic forms and Hecke theory and is also used by Templier (see (5.8) in \cite{templier}). Nevertheless, the difference consists merely of a character on $\mathfrak{o}_F^\times / \mathfrak{o}_F^{\times 2}$.}
The same consideration on the pretrace formula as in Section \ref{sec-amp} (without the amplifier) apply to this case as well.

Next, Hecke theory over number fields is best understood adelically. Given our classical motivation and for the sake of brevity, we prefer not to introduce the general formalism and additional notation, since it is not essential for the main argument of this paper. We refer to the detailed description of the classical and adelic theory given in Section 2 of \cite{shimura} in the related case of Hilbert modular forms. 

Instead, we note as in \cite{templier}, Section 5.5, that the action of a subalgebra of the Hecke algebra suffices for our purposes. Namely, to amplify the pretrace formula, we use only Hecke operators that are associated to principal ideals of $\mathfrak{o}_F$, whose action is explained for instance in Templier's article. By Chebotarev's density theorem (see Theorem 13.2 in \cite{neukirch}), these ideals make up a positive proportion of all prime ideals of $F$, it's numerical value depending only on $F$. As before, we assume that $\mathcal{O}$ is locally norm-maximal for simplicity and recall Remark \ref{rem-hecke}.

Define
\begin{displaymath} 
\mathcal{O} (m ; z, \delta) = \{ \gamma \in \mathcal{O}: |N_F(\nr( \gamma ))| = m, d (z, \varphi_0(\tilde{\gamma}) z) = O (\delta) \}.
\end{displaymath}
Using the notation of Section \ref{sec-amp} and following Sections 6.5 and 6.6 in \cite{templier}, we obtain a version of \eqref{eq-pt-ineq-split}. There is a certain sequence $y_m$ supported on positive integers less than $L^4$ such that
\begin{equation}\label{eq-amp-quat}
L^{2-\varepsilon} D^{-\varepsilon} |\phi(z)|^2 \ll_F S(\mu^\ast) \cdot \left( \sum_{m \ll L^4} \frac{y_m}{\sqrt{m}} \# \mathcal{O}(m; z, \delta) + 
(S(\mu^\ast) \delta)^{-\frac{1}{2}} \sum_{m \ll L^4} \frac{y_m}{\sqrt{m}} \# \mathcal{O}(m; z, \rho) \right).
\end{equation}
Furthermore, $y_m \ll 1$ and $\sum_m y_m \ll L^2$, as in (6.17) of \cite{templier}. 

As before, we are now faced with a counting problem. Note that we may count elements of 
$\mathcal{O} (m ; z, \delta)$ modulo units $\mathfrak{o}_F^\times$, since $F \subset \mathbb{R}$ is the centre of $A$. This remark is much more useful in the number field case than over the rationals, given that $\mathfrak{o}_F^\times$ is generally infinite.


\subsection{The counting argument}\label{sec-count-quat}

We now follow the argument in Section \ref{sect-disc-count}. We shall therefore show that the algebra generated by the elements we are counting is a proper algebra and thus a field by the tower rule.

Recall that by Lemma \ref{lemma-alg-in-vs}, which is independent of the ground field, it suffices to show that the $\mathbb{Q}$-\emph{vector space} spanned by $\bigcup_{1 \leq m \leq L^2 } \mathcal{O} (m ; z, 2 \delta)$ is proper. The following is an adapted version of Lemma \ref{lemma-proper}. By chance though, this lemma now gives bounds in terms of the parameter $\delta$ as well. We exploit the fact that $\operatorname{SO}(2)$ can only span a two dimensional vector space. This behaviour is not generic since orthogonal matrices in degree larger than $2$ can span the entire algebra of matrices over $\mathbb{R}$.

\begin{lemma}[Lemma \ref{lemma-proper} revisited] \label{lemma-proper-nf}
The $F$-vector space spanned by $\bigcup_{1 \leq m \leq L^{2}} \mathcal{O} (m ; z, 2 \delta)$ is proper, i.e.\ not equal to $A$, if $\delta \ll D^{1 - \varepsilon} L^{-6 - \varepsilon}$, where the implied constant depends only on $\varepsilon$.
\end{lemma}
\begin{proof}
Let $\gamma_1, \gamma_2, \gamma_3, \gamma_4 \in \bigcup_{1 \leq m \leq L^2 } \mathcal{O} (m ; z, 2 \delta)$. We wish to show that these elements are linearly dependent and we may assume without loss of generality that $\gamma_1 = 1$. Since the reduced trace gives a \emph{non-degenerate} bilinear form, it suffices to show that $s := \det( \tr(\gamma_i, \gamma_j) )_{i,j} = 0$.

By assumption, we have $\sigma_0(\nr(\gamma))^{-1/2} \cdot z^{-1} \varphi_0(\gamma_i) z = k_i + O(\delta)$, for some $k_i \in \operatorname{SO}(2)$, and thus
\begin{displaymath}
\frac{1}{\sigma_0 (\nr(\gamma_i) \nr(\gamma_j))^{1/2}} z^{-1} \varphi_0(\gamma_i \gamma_j) z = k_i k_j + O(\delta).
\end{displaymath}
Therefore,
\begin{displaymath}
\frac{1}{\prod_i \sigma_0(\nr(\gamma_i))} \sigma_0 \det(\tr(\gamma_i \gamma_j)_{i,j}) = \det(\tr (k_i k_j)_{i,j}) + O(\delta).
\end{displaymath}

Note that $\operatorname{SO}(2)$ spans only a $2$-dimensional vector space, which is easily seen using the standard parametrisation. It follows that $\det(\tr (k_i k_j)_{i,j}) = 0$, since we are considering the $4$-dimensional matrix space. Therefore, 
\begin{displaymath}
\sigma_0 \det(\tr(\gamma_i \gamma_j)_{i,j}) \ll \delta \prod_i |\sigma_0 \nr(\gamma_i)|,
\end{displaymath}
recalling that $\gamma_1 = 1$.

If $\sigma \neq \sigma_0$, then $\sigma \tr(\gamma_i \gamma_j) \leq 2 |\sigma \nr(\gamma_i \gamma_j)|^{1/2}$. Indeed, for the corresponding projections $\varphi_i, i\neq 0$, $\varphi_i(\xi)$ is an element of the real Hamilton quaternion algebra. For an arbitrary such element $a+ib+jc+kd$ in the usual notation with $a,b,c,d \in \mathbb{R}$, its trace is equal to $2a$ and its norm is $a^2 + b^2 + c^2 + d^2$, whence the inequality. Therefore, we deduce that
\begin{displaymath}
\sigma \det(\tr(\gamma_i \gamma_j)_{i,j}) \ll \delta \prod_i |\sigma \nr(\gamma_i)|.
\end{displaymath}

By the inequalities above, we have that
\begin{displaymath}
N_F(s) \ll \delta \prod_{i} |N_F(\gamma_i)| \ll \delta L^6.
\end{displaymath}
Recall that $\mathfrak{D}$ is the generator of the ideal in $\mathfrak{o}_F$ generated by $\{ \det\operatorname{tr}(x_i x_j) \mid x_i \in \mathcal{O}, i = 1,\ldots 4\}$. Since $\gamma_i \in \mathcal{O}$ for all $i$, it follows that $s = \mathfrak{D} \cdot x$ for some $x \in \mathfrak{o}_F$ and therefore $D \mid N_F(s)$. 

In conclusion, if $\delta \ll D^{-1 - \varepsilon} L^{6 - \varepsilon}$, then $s = 0$. By the non-degeneracy of the bilinear form given by the reduced trace, it follows that $\gamma_1, \ldots, \gamma_{p^2}$ are \emph{not} linearly independent.
\end{proof}

\begin{remark}
We remark that for $n \geq 3$ the $n \times n$ orthogonal matrices span the full space of real matrices. Thus, an application of the same proof as above in higher rank is bound to fail since the determinant $\det(\tr (k_i k_j)_{i,j})$ might be non-zero.
\end{remark}

Next, Lemma \ref{lemma-bound-in-field} goes through with the same proof if we can control the action of units. This cannot be done directly as in Lemma \ref{lemma-units-field} since $N_{E/F} (\xi) \in \mathfrak{o}_F^\times$ for any unit $\xi \in \mathcal{O}_E^\times$, and the group $\mathfrak{o}_F^\times$ is infinite for $F \neq \mathbb{Q}$. We can balance this out by recalling that we only need to count $\xi$ up to units in $\mathfrak{o}_F^\times$, since these act trivially on the upper half plane.

\begin{lemma}[Lemma \ref{lemma-units-field} revisited]\label{lemma-units-field-nf}
Let $E / F$ be an extension of degree $2$ that is a subfield of $A$ and let
$\mathcal{O}$ be an order of $A$. The number of units $\xi \in
\mathcal{O}^{\times} \cap E$ up to multiplication by units in $\mathfrak{o}_F$, such that $d (z, \varphi_0(\xi) z) \leq \delta$ for some $z \in \mathfrak{h}^2$, is $\ll_{F} (1+\delta)^2$.
\end{lemma}
\begin{proof}
We begin by investigating the quantity $(\tr_{E/F} \xi)^2 / N_{E/F} \xi$ and proving that it can only take finitely many values. For any embedding $\sigma \neq \sigma_0$, we have
\begin{displaymath}
	\sigma \left( \frac{(\tr_{E/F} \xi)^2 }{N_{E/F} \xi} \right) \in [0, 4],
\end{displaymath}
and the condition $d (z, \varphi_0(\xi) z) \leq \delta$ implies that
\begin{displaymath}
\sigma_0 \left( \frac{(\tr_{E/F} \xi)^2 }{N_{E/F} \xi} \right) \ll (1+\delta)^2,
\end{displaymath}
as we have seen already in the proof of Lemma \ref{lemma-proper-nf}.

Since $\xi \in \mathcal{O}_E^\times$, the maximal order, the quantity $(\tr_{E/F} \xi)^2 / N_{E/F} \xi$ must lie in $\mathfrak{o}_F$. Recall that the image of $\mathfrak{o}_F$ inside $\mathbb{R}^n$ under all embeddings is a discrete lattice. Since the image of $(\tr_{E/F} \xi)^2 / N_{E/F} \xi$ is bounded, it follows that the number of possibilities for the value of this quantity is bounded by $(1+\delta)^2$, up to a constant depending on $F$.

For the last step, recall Dirichlet's unit theorem, stating that $\mathfrak{o}_F^\times$ is a finitely generated group. This implies that $\mathfrak{o}_F^\times / (\mathfrak{o}_F^\times)^2$ is finite. Now if $\kappa \in \mathfrak{o}_F^\times$, then $N(\kappa \xi) = \kappa^2 N(\xi)$. Thus, if we are only counting $\xi \in \mathfrak{o}_F^\times \backslash \mathcal{O}_E^\times$, then the value of $N(\xi)$ can only lie in $\mathfrak{o}_F^\times / (\mathfrak{o}_F^\times)^2$. 

Nevertheless, we have
\begin{displaymath}
	\frac{(\tr_{E/F} (\kappa \xi))^2 }{N_{E/F} (\kappa \xi)} = \frac{(\tr_{E/F} \xi)^2 }{N_{E/F} \xi}.
\end{displaymath}
Since there are only finitely many possibilities for this quantity and finitely many for $N_{E/F}(\xi)$, if follows that there are only finitely many possibilities for $\tr_{E/F} \xi$. Finally, $\xi$ is determined up to the action of $\text{Gal}(E/F)$ (which has order 2) by its minimal polynomial. This polynomial is determined by the trace and norm of $\xi$. The lemma now follows by bookkeeping. 
\end{proof}

\begin{remark}
We remark that Lemma 6.4, part (i), of \cite{templier}, having the same ultimate goal as Lemma \ref{lemma-units-field-nf} in this paper, might not hold in general. Indeed, in the proof, the condition $u (z, \varphi_0(\xi) z) \leq \delta$ is said to be equivalent to \[\frac{\varphi_0(\xi)}{\sigma_0(N_A(\xi))^{1/2}}
\in zB (\delta) z^{- 1},\] for a single $\delta$-ball $B(\delta)$ around the identity. This cannot be true in general. For instance if $z = i$, corresponding to the identity matrix, and $\xi = \begin{psmallmatrix}
0 & -1 \\ 1 & 0
\end{psmallmatrix}$, then $u (z, \varphi_0(\xi) z) = 0$, yet $\xi$ lies far from the identity.

Our proof uses the idea of Lemma 6.3 in \cite{templier} and couples it with an application of Dirichlet's unit theorem, so as to apply the argument using characteristic polynomials, as in the higher degree case.
\end{remark}

Having adapted the lemmata above, we now complete the argument as in the Section \ref{sect-disc-count} and obtain the following counting result.

\begin{proposition}\label{prop-count-quat}
Let $\delta \ll D^{1 - \varepsilon} m^{-6 - \varepsilon}$ and let $z \in \mathfrak{h}^2$. Then $\# \mathcal{O} (m ; z, \delta) \ll m^{\varepsilon}$, where the implicit constant depends only on $\varepsilon$ and $F$. 
\end{proposition}

Compared to the counting results in Sections \ref{sect-disc-count} and \ref{sec-spec-count}, Proposition \ref{prop-count-quat} is conceptually more satisfying, since it seemingly blends together the spectral and discriminant aspect. Nevertheless, we still need a counting results for large $\delta$.

\begin{proposition}\label{prop-triv-count-quat}
We have the bound $\# \mathcal{O}(m; z, \rho) \ll m^{1 + \varepsilon}$.
\end{proposition}
\begin{proof}
The proof is essentially the same as the proof of Proposition \ref{prop-spec-triv-count}. Yet again, additional care must be taken when counting units at the end. For this, we apply the same argument as in Lemma \ref{lemma-proper-nf}, which simplifies since $L$ is now equal to $1$. In particular, if $\rho \ll D$, then the relevant units generate a proper subalgebra, which must be a field in this case. Since $D$ is an integer and $\rho$ is small enough, the condition is met. To conclude the proof, we apply Lemma \ref{lemma-units-field-nf}.
\end{proof}

\subsection{Proof of Theorem \ref{thm-main-nf}}

Considering the range of $m$ in the pretrace inequality \eqref{eq-amp-quat}, in order to use Proposition \ref{prop-count-quat}, we take $\delta \asymp D^{1 - \varepsilon} L^{-24 - \varepsilon}$. Then, applying Proposition \ref{prop-triv-count-quat}, we have
\begin{align*}
|\phi(z)|^2 &\ll_F D^{\varepsilon} L^{-2 + \varepsilon}  S(\mu^\ast) \left( \sum_{m \ll L^4} y_m m^{-1/2 + \varepsilon} + (S(\mu^\ast) D^{1 - \varepsilon} L^{-24 - \varepsilon})^{-\frac12} \sum_{m \ll L^4} y_m m^{1/2 + \varepsilon} \right) \\
&\ll S(\mu^\ast) L^{-2 + \varepsilon} D^\varepsilon \left( \left(\sum_m y_m^2 \right)^{1/2} \left( \sum_m \frac{1}{m} \right)^{1/2} + (S(\mu^\ast) D)^{-\frac12} L^{12} \cdot L^2 \sum_m y_m \right) \\
&\ll S(\mu^\ast) L^\varepsilon D^\varepsilon \left(L^{-1} + (S(\mu^\ast) D)^{-\frac12} L^{14} \right),
\end{align*}
where we make use of the properties of the sequence $y_m$. Optimising by setting the two terms in the last factor equal, we set $L = (S(\mu^\ast) D)^{\frac{1}{30}}$ and obtain a saving of $L^{-1}$. This implies that
\begin{displaymath}
|\phi(z)|^2 \ll_F S(\mu^\ast) \cdot (S(\mu^\ast) D)^{-\frac{1}{30} + \varepsilon}.
\end{displaymath} 
This proves Theorem \ref{thm-main-nf}.

\section{Remarks on the case of composite degree}\label{sec-composite}

It is not clear how to generalise the counting arguments, particularly those in the discriminant aspect (Section \ref{sect-disc-count}), to the case of algebras of general degree $n \in \mathbb{N}_{\geq 2}$. The main issue is the existence of non-commutative proper subalgebras in general. Yet the ideas in this paper suffice for proving the partial results given in Theorem \ref{thm-eichler} for special types of orders, which we sketch in this section.

More precisely, let $N$ be a positive integer and $\mathcal{O}_0(N)$ be an order of the division algebra $A$ over $\mathbb{Q}$ of degree $n$ such that, at all unramified primes $p$, its completion is of the form
\begin{displaymath}
\mathcal{O}_0(N)_p = \left\{ \gamma \in M_n(\mathbb{Z}_p) \mid \text{last row of } \gamma \equiv (0, \ldots, 0, \ast)  \operatorname{mod} N \mathbb{Z}_p \right\},
\end{displaymath}
up to conjugation. These orders are interesting from the point of view of newform theory. If $n = 2$, then these are precisely the Eichler orders, which have generally received much attention in the theory of automorphic forms. In higher degree, these orders form a (proper) subset of the orders that are the intersection of two maximal orders. It is important for the proof of Theorem \ref{thm-eichler} to note that the level of an order of type $\mathcal{O}_0(N)$ is $N^{n-1}$. 

Now let $n$ be odd and set $\mathcal{O} := \mathcal{O}_0(N)$ for simplicity. We can then make the same observation as in the proof of Lemma \ref{lemma-comm-delta}, bound \eqref{eq-comm-bound}. Namely, if $\gamma_1, \gamma_2 \in \bigcup_{1 \leq m \leq L} \mathcal{O}(m ; z, \delta)$, then
\begin{displaymath}
\nr(\gamma_1 \gamma_2 - \gamma_2 \gamma_1) \ll \delta L^2.
\end{displaymath}

The advantage of working with the family of $\mathcal{O}_0(N)$ is that $N \mid \nr(\gamma_1 \gamma_2 - \gamma_2 \gamma_1) $. Indeed, an easy computation shows that all commutators of $\mathcal{O}_0(N)_p$, where $p$ is unramified, have last row congruent to the zero vector modulo $N \mathbb{Z}_p$. The claim follows since the norm can be computed locally.

The remarks above imply that $N \ll \delta L^2$ or $\nr(\gamma_1 \gamma_2 - \gamma_2 \gamma_1) = 0$. Thus, if $\delta \ll N^{1 - \varepsilon} L^{-2 - \varepsilon}$, then the algebra generated by $\bigcup_{1 \leq m \leq L} \mathcal{O}(m ; z, \delta)$ is commutative.  Therefore, the same counting strategy employed in the rest of this article (counting in commutative fields) would give a strong bound in this case as well.

We apply the bounds to the amplifier as follows. Let $\delta \asymp N^{1 - \varepsilon} L^{-2 - \varepsilon}$. Then by \eqref{eq-pt-ineq-split} we have
\begin{displaymath}
|\mathcal{P}|^2 \cdot |\phi(z)|^2 \ll_n S(\mu^\ast) (LN)^\varepsilon \left(
|\mathcal{P}| + |\mathcal{P}|^2 L^{-(n-1)} + S(\mu^\ast)^{\frac{-1}{n(n-1)}} \cdot N^{-\frac{1}{2}} L |\mathcal{P}|^2 \frac{L^{(n-1)^2 n}}{L^{(n-1)n}}  \right).
\end{displaymath}
These are essentially the same computations as in Section \ref{sec-proof}. We deduce in the same way that
\begin{displaymath}
|\phi(z)|^2 \ll S(\mu^\ast) (LN)^\varepsilon (L^{-1} + S(\mu^\ast)^{\frac{-1}{n(n-1)}} \cdot N^{-\frac{1}{2}} L^{1 + n(n-1)(n-2)}).
\end{displaymath}
We may simplify one of the exponents of $L$ by noting that $1 + n(n-1)(n-2) \leq n^3 - 1$ (for $n \geq 2$). We then find an optimal value of $L \asymp S(\mu^\ast)^{1/n^4(n-1)} N^{1/2n^3}$, so that
\begin{displaymath}
|\phi(z)|^2 \ll S(\mu^\ast)^{1 - \frac{1}{n^4(n-1)} + \varepsilon} \cdot N^{-\frac{1}{2n^3} + \varepsilon}.
\end{displaymath}
This proves Theorem \ref{thm-eichler}.

The sup-norm bound we obtain is unfortunately not uniform in the full volume aspect, since we did not include the discriminant of $A$ in the bounds. It is possible to include ramified primes $p$ where $A_p$ is a division algebra, since then $p$ divides the norm of commutators. Nevertheless, in the composite degree case considered here there is also the possibility of $A_p$ being a more general matrix algebra over a division algebra, which we are not able to treat using this approach. It would certainly be interesting to at least extend these bounds to include the full discriminant, but even more so to find a more flexible argument to treat arbitrary orders. 

For example, the argument does not apply to the larger family of generalised Eichler orders, which we define to be intersections of two maximal orders. One example in degree 4, which is also a hereditary order, has the form
\begin{displaymath}
\mathcal{O}_p =
\begin{pmatrix}
\mathbb{Z}_p & \mathbb{Z}_p & \mathbb{Z}_p & \mathbb{Z}_p \\
\mathbb{Z}_p & \mathbb{Z}_p & \mathbb{Z}_p & \mathbb{Z}_p \\
p\mathbb{Z}_p & p\mathbb{Z}_p & \mathbb{Z}_p & \mathbb{Z}_p \\
p\mathbb{Z}_p & p\mathbb{Z}_p & \mathbb{Z}_p & \mathbb{Z}_p \\
\end{pmatrix},
\end{displaymath}
at a prime $p$. In this case, the norm of a commutator need not be divisible by $p$.

\printbibliography

\end{document}